\documentclass[12pt]{amsart} 
\usepackage{amsmath,latexsym,amssymb}
\begin{document}
\newtheorem{theorem}{Theorem}[section]
\newtheorem{definition}[theorem]{Definition}
\newtheorem{proposition}[theorem]{Proposition}
\newtheorem{lemma}[theorem]{Lemma}
\newtheorem{remark}[theorem]{Remark}
\newtheorem{corollary}[theorem]{Corollary}
\newtheorem{question}[theorem]{Question}
\newtheorem{example}[theorem]{Example}
\newtheorem{notation}[theorem]{Notation}
\newtheorem{claim}[theorem]{Claim}
\def\cl{\begin{claim}\upshape}
\def\ecl{\end{claim}}
\def\prcl{\par\noindent{\em Proof of Claim: }}
\def\bem{\begin{remark}\upshape}
\def\ebem{\end{remark}}
\def\nota{\begin{notation}\upshape}
\def\enota{\end{notation}}
\def\defn{\begin{definition}\upshape}
\def\edefn{\end{definition}}
\def\thm{\begin{theorem}}
\def\ethm{\end{theorem}}
\def\lmm{\begin{lemma}}
\def\elmm{\end{lemma}}
\def\qed{\hfill$\quad\Box$}
\def\pr{\par\noindent{\em Proof: }}
\def\prbeh{\par\noindent{\em Proof of Claim: }}
\def\cor{\begin{corollary}}
\def\ecor{\end{corollary}}
\def\frag{\begin{question}\upshape}
\def\efrag{\end{question}}
\def\prop{\begin{proposition}}
\def\eprop{\end{proposition}}
\def\beh{\par\noindent{\bf Claim: }\em}
\def\ebeh{\par\noindent\upshape}
\def\bsp{\begin{example}}
\def\ebsp{\end{example}}
\def\a{{\sigma}}
\def\si{{\sigma}}
\def\F{{ \mathbb F}}
\def\N{{\mathbb N}}
\def\Z{{\mathbb Z}}
\def\Q{{\mathbb Q}}
\def\L{{\mathcal L}}
\def\P{{\mathcal P}}
\def\I{{\mathcal I}}
\def\B{{\mathcal B}}
\def\C{{\mathcal C}}
\def\IC{{\mathbb C}}
\def\R{{\mathcal R}}
\def\IR{{\mathbb R}}
\def\M{{\mathcal M}}
\def\I{{\mathcal I}}
\def\D{{\mathcal D}}
\def\O{{\mathcal O}}
\def\Ma{{\mathfrak m}}
\def\U{{\Upsilon}}
\def\H{\mathcal{H}}
\title{B\'ezout domains and lattice-valued modules.}
\date{\today}
\author{Sonia L'Innocente}
\address{\hskip-\parindent 
L'Innocente Sonia\\
School of Science and Technology - Department of Mathematics\\
University of Camerino\\
Via Madonna delle Carceri 9, 62032 Camerino (MC), Italy.}
\email {sonia.linnocente@unicam.it}
\author{Fran\c coise Point$^{\dagger}$}
\address{\hskip-\parindent
Fran\c coise Point\\
Department of Mathematics (Le Pentagone)\\
University of Mons\\
20, place du Parc, B-7000 Mons, Belgium.}
\email {point@math.univ-paris-diderot.fr}
\thanks{Both authors gratefully acknowledge the support of the project FIRB2010 (Italy), of which the first author is the principal investigator.\\
${}^{\dagger}$Research Director at the \lq\lq Fonds de la Recherche
Scientifique FNRS-FRS\rq\rq}

\begin{abstract}
Let $B$ be a commutative B\'ezout domain $B$ and let $MSpec(B)$ be the maximal spectrum of $B$. We obtain a Feferman-Vaught type theorem for the class $Mod_{B}$ of $B$-modules. We analyse the definable sets in terms, on one hand, of the definable sets in the classes $Mod_{B_{\M}}$, where $B_{\M}$ ranges over the localizations of $B$ at $\M\in MSpec(B)$, and on the other hand, of the constructible subsets of $MSpec(B)$. When $B$ has good factorization, it allows us to derive  decidability results for the class $Mod_{B}$, in particular when $B$ is the ring $\widetilde \Z$ of algebraic integers or the one of rings $\widetilde \Z\cap \IR,\;\widetilde \Z\cap \Q_p$.
\end{abstract}

\maketitle
\section{Introduction} 
\par Let $B$ be a commutative B\'ezout domain with $1$ and let $B^{\star}:=B\setminus \{0\}$. Let $MSpec(B)$ be the space of maximal ideals of $B$ endowed with the Zariski topology. (Basic closed sets are $V(a):=\{\M\in MSpec(R): a\in \M\}$ and constructible subsets of $MSpec(B)$ are the elements of the Boolean algebra generated by the basic closed subsets.) 
\par In the class $Mod_{B}$ of all $B$-modules $M$, we will describe the definable subsets of $M$ in terms of the definable sets  in each localization $M_{\M}$, $\M\in MSpec(B)$,
 and of constructible subsets of $MSpec(B)$. This description of definable subsets can be seen as a Feferman-Vaught type result and this is the content of our main Theorem \ref{ppelB}.  
 We will work in a definable expansion of the language of $B$-modules, adding to the usual language of modules unary predicates for submodules indexed by the group of divisibility $\Gamma(B)$ of $B$. 
\par A key intermediate step in the proof of Theorem \ref{ppelB} is a positive quantifier elimination result in that expansion of the language (showing that any positive primitive formula is equivalent to a conjunction of atomic formulas) in the class $Mod_B$ when $B$ is a valuation domain (Theorem \ref{ppel}). Even though, this result was essentially known, in order to apply it when $B$ is a B\'ezout domain, we need to control what happens when we localize $B$ at maximal ideals.
\par For $a, b\in B^{\star}$, recall that the Jacobson radical relation $a\in rad(b)$ holds if $a$ belongs to every maximal ideal that contains $b$ and  denote by $gcd(a,b)$ the greatest common divisor of $a,\,b$. 
\par In \cite[2.10]{DM}, L. van den Dries and A. Macintyre defined the property for a B\'ezout domain $B$ to have {\it good factorisation}, namely if for all $a,\;b\in B^{\star}$, there are $c,\;a_{1}\in B$ such that $a=c.a_{1}$ with $gcd(c,b)=1$ and $b\in rad(a_{1})$. Then they observed that in case $B$ has good factorization, the constructible subsets of $MSpec(B)$ are either basic closed sets or basic open sets \cite[Lemma 2.12]{DM}. Furthermore, they show that one can encode the properties of the Boolean algebra of constructible subsets of $MSpec(B)$ using the Jacobson radical relation.
\par Then, under the further assumption that $B$ has good factorization, we consider questions of decidability for the class $Mod_{B}$.
First we make explicit in our context the property for a countable ring $B$ to be {\it effectively given} \cite[Chapter 17]{Pr}; we call the resulting assumption on $B$, {\it assumption (EF)} (Definition \ref{EF}).
\par We show that the class $Mod_{B}$ is decidable (Proposition \ref{dec}), under the following hypotheses on $B$: $B$ is countable and it satisfies assumption (EF), it has good factorization, each quotient $B/\M$ is infinite, for $\M\in MSpec(B)$, and the Jacobson radical relation is recursive.  
\par B\'ezout domains with good factorization include the class of good Rumely domains introduced by L. van den Dries and A. Macintyre \cite{DM}. Even though, their work takes place in the context of rings, there are similarities in their approach and the one we are taking. They show that the class of good Rumely domains admits quantifier elimination in the language of rings extended by a family of radical relations \cite[Theorem 3.14]{DM}.  
There are three main ingredients in their proof: Rumely local-global principle, the quantifier elimination result for non-trivial valuation rings with algebraically closed fraction fields (following from a theorem of A. Robinson) and the fact that the constructible subsets of the maximal spectrum of a such ring form an atomless Boolean algebra.
They axiomatize the class of good Rumely domains and retrieve the former result of van den Dries \cite{vdD}  on the decidability of the ring of algebraic integers.
\medskip
\par The plan of the paper is as follows. In section 2, we recall the basic notions of the model theory of modules (or abelian structures) that we will use and 
 the properties of the group of divisibility $\Gamma(B)$ of a B\'ezout domain. 
\par In section 3, for $A$ a valuation domain, we revisit a quantifier elimination result in the class $Mod_A$, adding to the module language unary predicates for certain pp definable submodules, following the approach of B\'elair--Point \cite[Proposition 4.1]{BP}. This result was essentially known but we need the additional property that for $A$ of the form $B_\M$, where $B_{\M}$ is the localization of a B\'ezout domain $B$ at a maximal ideal, to any pp $\L_{B_{\M}}$-formula one can associate a constructible subset of $MSpec(B)$ over which the elimination is uniform.  
\par In section 4, we give a direct proof of a decidability result in the case of valuation domains with infinite residue field. This was previously done in \cite{PPT} in the case of an archimedean value group and later extended in \cite{Gr} in the general case (and furthermore without the assumption on the residue field).
The proof given here is more algebraic than the one in \cite{PPT} which is of a more geometrical nature. We apply the result in section 5 in the case of B\'ezout domains.
\par  In section 5, we first relate the property of having good factorization for $B$ to properties of $MSpec(B)$ and note that such ring is {\it adequate}, a better known property. Then we prove our main theorem, a Feferman-Vaught type result for the class $Mod_{B}$, which takes a simpler form in the case where $B$ has good factorization. 
\par We derive a decidability result when $B$ is an effectively given countable B\'ezout domain where the Jacobson radical relation is recursive, assuming that the quotient $B/\M$ is infinite for any maximal ideal $\M$. Note that L. Gregory observed that if the theory of all $B$-modules is decidable, then the prime radical relation is recursive \cite[Lemma 3.2]{Gr}. For good Rumely domains, the prime radical relation and the Jacobson radical relation $rad$ always coincide \cite{PS} and so our hypothesis on the Jacobson radical relation is justified in view of Gregory's result.
\par In section 6, we apply our decidability result to the case where $B$ is an effectively given good Rumely domain e.g. $\widetilde \Z$ and to the cases $\widetilde \Z\cap \IR,\;\widetilde \Z\cap
\Q_{p}$ \cite{PS91}. 
\par In the last subsection, we discuss the case when $B$ is either the ring of holomorphic functions over $\IC$ or the integral closure of that ring. Of course in this case, the ring is uncountable (and so
the language of modules is uncountable), but also it is not known to satisfy the Rumely local-global principle \cite[5.6]{DM}. However since these rings has good factorization, we still have a manageable description of definable subsets in that class, using Corollary \ref{cons}.
\par Then, in the last section, we introduce the notion of $\ell$-valued $B$-modules, in view of future work. When $B$ is a valuation domain, we get back the more classical notion of valued modules \cite{Co}, \cite{F}.
\smallskip
\par Finally let us note that since this paper was submitted other works on decidability of the theory of modules over B\'ezout domains appeared; see for instance {\it Decidability of the theory of modules over B\'ezout domains with infinite residue fields, arXiv:1706:08940} by L. Gregory, S. L'Innocente, G. Puninski and C. Toffalori. 
\section{Preliminaries}\label{prelim}
\par Throughout the paper, all our rings $B$ will be commutative integral domains with $1$. Let $B^{\star}:=B\setminus \{0\}$. Then $B$ is B\'ezout if every finitely generated ideal is principal, equivalently if 
for any $a, b\in B^{\star}$, there exist $c, u, v, a_{1},a_{2}\in B$, such that $c=a.u+b.v$ and $a=c.a_{1}$, $b=c.b_{1}$ (the B\'ezout relations). We set $c=gcd(a,b)$ and we denote by $(a:b):=a_{1}$, a generator of the ideal $(a):(b):=\{u\in B:\;b.u\in (a)\}$; note that these two elements $c,\;a_{1}$ of $B$ are defined up to an invertible element. We can also define the least common multiple of two elements $a,\;b$, denoted by $lcm(a,b)$. It is easily checked that $a.b=gcd(a,b).lcm(a,b).u$, where $u$ is an invertible element of $B$.
\par In the next two subsections, we will quickly review some basic facts on one hand  the group of divisibility of a B\'ezout domain and on the other hand abelian structures.
\subsection{Group of divisibility} $\;$
\par Let $B$ be a B\'ezout domain and denote by $Q(B)^{\times}$ the multiplicative group of the field of fractions of $B$ and by $U$ the subgroup of units (equivalently of invertible elements) of  $B$. 
\par Recall that the group of divisibility $\Gamma(B)$ of $B$ is the quotient of $Q(B)^{\times}$ by $U$, more formally  $\Gamma(B)=(Q(B)^{\times}/U,\cdot ,1)$, where $1$ denotes the neutral element and is the coset $1\cdot U$. Denote by $\Gamma^+(B):=\{a.U:\;a\in B^{\star}\}$. 
One can define a partial order on $\Gamma(B)$ by: $a.U\leq b.U$ iff $b.a^{-1}\in B$. Endowed with this partial order, $\Gamma(B)$ becomes a partially ordered group \cite[1.2]{Gl}. Note that $\Gamma^{+}(B)$ corresponds to the submonoid of positive elements of $\Gamma(B)$ (namely those bigger than or equal to $1$). Recall that since $B$ is B\'ezout, this order is a lattice order and so $\Gamma(B)$ is an abelian lattice-ordered group (in short, $\ell$-group) \cite[Chapter 3, Proposition 4.5]{FS}. 
For $a,\; b\in B^{\star}$, we can explicitly define the lattice operations on the set of positive elements $\Gamma^+(B)$ as follows: $a.U\wedge b.U:=gcd(a,b).U$ and $a.U\vee b.U=lcm(a,b).U$. One then extends the lattice operations to the elements of $Q(B)^{\times}$.
\par The map sending $a\in B^{\star}$ to the coset $a.U$ is called an $\ell$-valuation, generalizing the notion of a valuation map on an valuation domain.
\defn {\rm \cite[Definition 1]{RY}} \label{Rump} Let $D$ be an integral domain, $(\Gamma,\wedge,\cdot,1)$ be an $\ell$-group and set $\bar \Gamma:=\Gamma\cup\{\infty\}$ with $\infty\cdot a=\infty\geq a$ for $a\in \bar \Gamma$.
Then the map $v:\;D\rightarrow \bar \Gamma$ is an $\ell$-valuation if it satisfies the
properties $(1)$ up to $(3)$ below:
for all $a,\;b\in D$,
\begin{enumerate}
\item $v(a+b)\geq v(a)\wedge v(b)$,
\item $v(a.b)=v(a).v(b)$,
\item $v(1)=1$ and $v(0)=\infty$.
\end{enumerate}
\edefn
\par One can show that every $\ell$-valuation on an integral domain $D$ has a unique extension to its field of fractions $Q(D)$ \cite[Corollary to Proposition 2]{RY}.
\medskip
\par When $B$ is a B\'ezout domain, the map $v:B^{\star}\rightarrow\Gamma(B): a\mapsto a.U$ is an $\ell$-valuation on $B$.
For $\M\in MSpec(B)$, the space of maximal ideals of $B$, we denote by $B_{\M}$ the localization of $B$ at $\M$. Let $U_{\M}$ denote the subgroup of invertible elements of $(B_{\M}^{\star},\cdot,1)$.
Since $B$ is a B\'ezout domain, $B_{\M}$ is a valuation domain and $\Gamma(B_{\M})$ is an abelian totally ordered group. 
\par It is well-known that any $\ell$-group is isomorphic to a subdirect product of totally-ordered abelian groups \cite[Lemmas 3.2.5, 3.5.4]{Gl}. In our framework, it is useful for us to describe that isomorphism as follows. We will indicate that we consider a subdirect product by using $\prod^s$.
\lmm \label{emb} \label{embgr}  The map $f:\;a.U\rightarrow (a.U_{\M})_{\M\in MSpec(B)},$ with $a\in B$, 
induces an isomorphism between the lattice-ordered monoid $(\Gamma^+(B),.,\wedge,1)$ and \\ 
$\prod_{\M\in MSpec(B)}^s (\Gamma(B_{\M})^+,.,\wedge, 1)$, where $(\Gamma^+(B_{\M}),.,\wedge, 1)$, $\M\in MSpec(B)$ are totally ordered monoids. 
It can be extended to an isomorphism of lattice-ordered
groups, that we will still denote by $f$, from $(\Gamma(B),.,\wedge,1)$ to $\prod_{\M\in MSpec(B)}^s (\Gamma(B_{\M}),.,\wedge, 1)$.
\qed
\elmm
\par In particular, we denote by $f_{\M}:\Gamma(B)^{+}\rightarrow\Gamma(B_{\M})^{+}:\;a.U\rightarrow a.U_{\M}$ with $a\in B^{\star}$. This induces an $\ell$-valuation $v_{\M}$ on $B_{\M}$, with values in $\Gamma(B_{\M})^+$ as follows: let $a\in B$, $s\in B\setminus \M$, then  $v_{\M}(a.s^{-1}):=a.s^{-1}.U_{\M}=a.U_{\M}=f_{\M}(a.U)$.
\subsection{Abelian structures} \label{prelim2} $\;$
\par We will consider the class of (right) $B$-modules endowed with a family of subgroups. Let $\L_{B}:=\{+,-,0,\cdot a; a\in B\}$ be the language of $B$-modules, where $\cdot a$ denotes scalar multiplication by $a$, and let $\L_{B,V}$ be the language of $B$-modules expanded  by a set $\{V_{\delta}; \delta\in \Delta\}$ of unary predicates (namely unary relation symbols); i.e. $\L_{B,V}:=\L_{B}\cup\{V_{\delta}; \delta\in \Delta\}$ (and $\Delta$ is some index set). (When the ring $B$ is clear from the context, we will drop the subscript $B$).
\par Given a $B$-module $M$ , we will consider its expansion $M_{V}$ by a prescribed family of subgroups $V_{\delta}(M)$, $\delta\in \Delta$; this is an instance of an abelian structure \cite{Fis}, \cite[Chapter 3, 3A]{Pr}. E. Fisher in his thesis had extended to the class of abelian structures most of the classical results on the model theory of theory of modules. At the beginning of \cite{Z}, M. Ziegler pointed out that most results in the model theory of modules still hold in this larger setting. Later, in \cite[section 1.9]{KP}, T. Kucer\u{a} and M. Prest described a way to view any abelian structure as a module over a certain path algebra associated with the language. This point of view has the advantage of staying in the classical framework of modules but the disadvantage of changing the ring. 
\par We are interested in describing the definable subsets of such $\L_{V}$-structure $M_{V}$ and we will use the Baur--Monk pp elimination theorem, namely any $\L_{V}$-formula is equivalent to a boolean combination of pp formulas and invariant sentences \cite[Appendix A.1]{Hod}. Moreover, such pp elimination is effective and uniform in the class of $\L_{V}$-structures.
\par Let $p^{+}(x)$ (respectively $p^{-}(x)$) be a set of pp formulas (in one variable, in $\L_{V}(A)$, where $A$ is a set of parameters) such that any pp formula either belong to $p^{+}$ or to $p^{-}(x)$ (but not to both) and such that given any finite subset $E$ of formulas in $p^{+}(x)$ there is a module $M_{V}$ (containing $A$) and an element $m\in M_{V}$ such that $\phi(m)$ holds in $M_{V}$ for any $\phi\in E$. We will call $p^+$ a pp type. Denote by $\neg p^{-}(x)$ the collection of negations of formulas in $p^{-}(x)$. Then a type $p(x)$ (in one variable) is of the form $p^{+}\cup \neg p^{-}(x)$ and we say it is realized in a structure $M_{V}$ if there is an element $m$ in $M_{V}$ such that $\phi(m)$ (respectively $\neg \phi(m)$) holds, for any $\phi\in p^{+}$ (respectively for any $\phi\in p^{-}$). 
\par A pure-injective abelian structure $M_{V}$ is an abelian structure where every pp type is realized (with parameters in $A$ with $\vert A\vert\leq \vert \L_{V}\vert$ \cite[Theorem 3.1]{Z}); $M_V$ is indecomposable when one cannot decompose it as $M_{1}\oplus M_{2}$, with $M_{1},\,M_{2}$ non-zero.
With a type $p$ (over the empty set), one may associate a unique minimal pure-injective structure $H(p)=H(m)$ \cite[Theorem 3.6]{Z} with $m\in H(p)$ such that $p^{+}(x)$ is exactly the set of pp-formulas satisfied by $m$ in $H(p)$. One says that $p$ is indecomposable if $H(p)$ is indecomposable.
\par A basic result in the model theory of modules which has been adapted to the setting of abelian structure is the following: any abelian structure $M_{V}$ is elementarily equivalent to a direct sum of pure-injective indecomposable abelian structures \cite[Corollary 6.9]{Z}. This result (for classical $B$-modules) has led M. Ziegler to associate with the class of $B$-modules, a topological space $Zg_{B}$ (the Ziegler spectrum) whose points are isomorphism types of non-zero pure-injective indecomposable modules and basic open sets (denoted by $[\phi/\psi]$, where $\phi,\;\psi$ are pp formulas in one variable with $\psi\rightarrow \phi$) consist of the (isomorphism types of) pure-injective indecomposable modules $M_V$ where the index of $\psi(M_V)$ in $\phi(M_V)$ is strictly bigger that $1$ \cite[Corollary 6.13]{Z}. This space is quasi-compact.
\medskip
\par Recall that two pp-formulas $\phi,\;\psi$ are said to be equivalent if in all modules $M$, $\phi(M)=\psi(M)$. Denote by $L(B)$ the set of pp $\L$-formulas quotiented by this equivalence relation. One can show that $L(B)$ forms a lattice \cite[Section 2]{PT} and in case $B$ is B\'ezout, this lattice is distributive  \cite[Fact 2.4]{PT}.
On that lattice of pp formulas, we will consider the duality functor $D$ which transforms a formula in the category of left $B$-modules to a formula in the category of right $B$-modules \cite[chapter 8]{Pr}. Since $B$ is commutative, it will allows us to simply the description of $Zg_{B}$.
\medskip
\par It has been observed that the group $\Gamma(B)$ of divisibility of $B$ reflects (some of) the model-theoretic properties of the class of $B$-modules. For instance in \cite[Theorem 7.1]{PT}, they show that $L(B)$ has no width if and only if $\Gamma(B)$ contains a densely ordered subchain. 
\par When $B,\;B'$ are two valuation domains (and so $\Gamma(B),\;\Gamma(B')$ are totally ordered abelian groups, also called the value groups), L. Gregory observed that the Ziegler spectrum $Zg_{B},\;Zg_{B'}$ are homeomorphic if and only if the value groups $\Gamma(B),\;\Gamma(B')$ are isomorphic \cite[Corollary 3.3]{Gr}.
\subsection{Pr\"ufer domains}\label{prufer}
\par Recall that a positive primitive (pp) $\L_V$-formula $\phi(\bar x)$, $\bar x:=(x_{1},\ldots, x_{n})$, $n\geq 1$, is an existential formula of the form: $\exists \bar y\; (\bar y.A=\bar x. C\;\&\;\bigwedge_{\delta\in \Delta} V_{\delta}(\bar y\bar x.C_{\delta}))$, where $\bar y=(y_{1},\ldots, y_{m})$, $m\geq 1$, $A$ is a $m\times k$-matrix, $C$ a $n\times k$ matrix, $k\geq 1$, $C_{\delta}$ a $(n+m)\times 1$ matrix, all with coefficients in $B$ and $\Delta$ a finite subset of $\Gamma$. By $\phi(M)$ we denote the submodule of $M^n$ consisting of the tuples $\bar u$ of elements of $M^{n}$ such that  for some $\bar w\in M^{m}$ we have $\bar w.A=\bar u. C\;\&\;V_{\delta}(\bar w\bar u.C_{\delta})$. Let $\phi(x),\;\psi(x)$ be two pp formulas in one free variable, then an invariant sentence is a sentence of the form $(\phi/\psi)> n$ that expresses that the index of the subgroup $\psi(M)\cap \phi(M)$ in the subgroup $\phi(M)$ is strictly greater than $n$, $n\in \N$. 
\medskip
\par Let us recall a former result of G. Puninski\u{i} \cite{Pun} on the special equivalent form of  pp formulas over a B\'ezout domain $B$ (or more generally over Pr\"ufer domains). 
\par Any pp $\L_{B}$-formula $\phi(\bar x)$ is equivalent  
to the pp $\L_{B}$-formula:
$\exists \bar y\; (\bar y.S=\bar x\;\&\;\bar y.\bar r=0)$, where $\bar y=(y_{1},\ldots, y_{m})$, $S$ is a $m\times n$-matrix, $\bar r$ a $m\times 1$ matrix,  all with coefficients in $B$. 
\par In view of this result, given a $B$-module $M$, $\delta\in \Gamma(B)^+$ and $a\in B$ with $v(a)=\delta$, we define $V_\delta(M):=\{m\in M:\;\exists y\in M\;m=y.a\}$. 
This is well-defined since for $a,\;b\in B^*$ with $v(a)=v(b)$, we have $a.b^{-1}\in B$. Since $B$ is commutative, $V_{\delta}(M)$ is not only a pp definable subgroup but an $B$-submodule. Finally note that the formula $V_{v(a)}(m.b)$ is either equivalent to the formula $m=m$ if $v(b)\geq v(a)$ or if $v(b)<v(a)$ to the formula $\exists y_1\exists y_2\exists z\;(m=y_1+y_2\;\&\;y_1.b=0\;\&\; y_2=z.a.b^{-1})$. (Indeed, in case $v(b)<v(a)$ and $V_{v(a)}(m.b)$ holds, there exists $n\in M$ such that $m.b=n.a=n.a.b^{-1}.b$ and so $m=m-n.a.b^{-1}+n.a.b^{-1}$ with $(m-n.a.b^{-1}).b=0$.
In the formalism used in \cite{PPT}, we will abbreviate $\exists y\;y.(a.b^{-1})=x$ by $a.b^{-1}\vert x$ and we will write the formula $V_{v(a)}(m.b)$ as $(b.x=0)+(a.b^{-1}\vert x)$, or equivalently $(b.x=0)+V_{v(a).v(b)^{-1}}(x)$.
\subsection{An example: the ring of holomorphic functions}$\;$\label{hol}
Let $\H(\IC)$ be the ring of holomorphic functions over $\IC$. A key ingredient is the Weierstrass factorization theorem \cite[Theorem 15.10]{Rudin}. One defines the functions $E_{0}(z)=1-z,\;E_{p}(z)=(1-z).exp\{z+z^{2}+\cdots+\frac{z^{p}}{p}\}$. 
Letting $(z_n)_{n\in \omega}$ be a sequence of elements of $\IC\setminus\{0\}$ (possibly with repetitions) such that $\vert z_n\vert\rightarrow \infty$, the infinite product $P(z):=\prod_{n=1}^{\infty} E_{n-1}(\frac{z}{z_n})$ belongs to $\H(\IC)$, the zeroes of $P$ are exactly the $z_n$'s and if $z_n$ occurs $m$ times in $P(z)$, then $z_n$ is a zero of $P(z)$ of multiplicity $m$ \cite[Theorem 15.9]{Rudin}. 
\bsp \label{hol} The ring $\H(\IC)$ is a B\'ezout domain \cite[Theorem 1]{H}, \cite[Theorem 15.15]{Rudin} (and we indicate a proof in  subsection \ref{entire}).  Let us describe its group of divisibility.
Given any element $f\in \H(\IC)$, we define the multiplicity function $\mu_{f}:\IC\rightarrow \N$ sending $z\in \IC$ to the multiplicity of $z$ as a zero of $f$.
Set $\daleth^+:=\{\mu_f:\;f\in \H(\IC)\}.$ We have that $\mu_{f.g}=\mu_{f}+\mu_{g}$; 
$\daleth^+$ forms a commutative monoid (w.r.to $+$) and can be endowed with a partial order: $\mu\leq \nu$ iff $\forall z\;\mu(z)\leq \nu(z)$, for $\mu,\;\nu\in \daleth^+$. 
This partial order reflects the divisibility relation in $\H(\IC)$: $f\vert g$ in $\H(\IC)$ iff $\mu_{f}\leq \mu_{g}$ and $f$ is invertible iff $\mu_{f}=0$. 
Denote by $\daleth$ the group generated by $\daleth^+$; it is easy to see that $(\daleth,\leq)$ is an $\ell$-group and that it is isomorphic to $\Gamma(\H(\IC)).$
\par Using the above description of the group of divisibility of $\H(\IC)$ and \cite[Theorem 7.1]{PT}, one can easily see that the lattice of pp formulas over $\H(\IC)$ has no width  \cite[Example 6.3]{PT}.
Indeed, choose two elements $f,\;g\in \H(\IC)$ with the same infinite (discrete) subset of zeroes: $Z(f)=Z(g)=\{z_{n}:\;n\in \N^{\star}\}$ and such that $\mu_g(z_n)=2n.\mu_f(z_n)$.
Then using Weierstrass factorization theorem recalled above, 
there exists $h\in \H(\IC)$ such that $\mu_{h}(z_{n}):=\mu_f(z_n)+ \lfloor \frac{\mu_{g}(z_n)-\mu_{f}(z_n)}{2}\rfloor$. Then $\mu_f<\mu_h<\mu_g$ and the strict inequality holds simply because $\lim_{n\rightarrow \infty} \mu_g(z_n)-\mu_f(z_n)=+\infty$, and so we may re-apply the same procedure to both: $(\mu_g,\mu_h),\;(\mu_h,\mu_f)$.
\ebsp

\section{Quantifier elimination for valuation domains}
\par Let $B$ be a B\'ezout domain and let $\Gamma(B)$ be its group of divisibility. Recall that the space $MSpec(B)$ of maximal ideals of $B$ is endowed with the Zariski topology where a basis of closed subsets is given by $V(a):=\{\M\in MSpec(B):\;a\in \M\}$. The {\it constructible} subsets of $MSpec(B)$ are the elements of the Boolean algebra generated by the basic closed subsets.
\par We have the following relationship between the lattice generated by these basic closed subsets of $MSpec(B)$ and the lattice of principal ideals of $B$. Let $a, b\in B^{\star}$. Then $V(a)\cap V(b)=V(gcd(a,b))$, $V(a)\cup V(b)=V(a.b)$, $V(a)\setminus V(b)\subset V((a:b))$.
\medskip
\par Let $\L_V$ be the language of $B$-modules expanded with a set of unary predicates $V_\delta$ indexed by the submonoid $\Gamma^+$ of positive elements of the group $\Gamma:=\Gamma(B)$. 
Let $M$ be an $B$-module and for $\delta=v(a)$, set $V_\delta(M):=\{m\in M:\;\exists n\in M\;m=n.a\}$, where $a\in B^{\star}$. 
\medskip
\par In the Theorem below, we consider a pp $\L_V$-formula $\exists
\bold x \,\phi(\bold x,\bold z)$ in the language $\L_{V}$, where $\phi(\bold x,\bold z)$ is a conjunction
of atomic formulas and we will show that this formula is equivalent in any $B_{\M}$-module to a
conjunction of atomic formulas. In case of the pure module language, that property is called {\it positive quantifier elimination} (or {\it elim-$Q^+$} following the terminology of \cite[page 319]{Pr}) and implies structural properties on indecomposable pure-injective modules \cite[Corollary 16.7]{Pr}.  
We could have applied Puninski\u{i}'s result on Pr\"ufer domains (recalled in subsection \ref{prelim2}) and the fact that over a valuation domain a matrix is conjugated to a diagonal matrix. 
However since we are ultimately interested in describing definable subsets in the class of $B$-modules (see section \ref{FV}), we will consider classes of modules over valuation domains $B_{\M}$, where $\M$ varies over $MSpec(B)$. So we will use that for a given pp $\L_{B_{\M}}$-formula, the elimination is uniform on a certain constructible subset of $MSpec(B)$. This is why we chose to give a self-contained and direct proof of that result.
\par We will use the isomorphism between $(\Gamma(B)^+,.,\wedge,1)$ and the subdirect product $\prod_{\M\in MSpec(B)}^s (\Gamma(B_{\M})^+,.,\wedge, 1)$ (Lemma \ref{embgr}). Let us denote the image of $\delta\in \Gamma(B)$ in $\Gamma(B_{\M})$ by $\delta_{\M}$.
Note that if $\delta_{\M}:=v_{\M}(a)$, $a\in B^{\star}$, and $M$ is a $B_\M$-module, we get that $V_{\delta}(M)=V_{\delta_{\M}}(M)$.
\bem \label{rad} {\rm \cite[section 1]{vdD}}
In the localization $B_{\M}$, we have that $a\vert b$ iff $v(a)\leq v(b)$ iff $(a:b)\notin \M$ iff $\M\notin V((a:b))$.
\par Since $B_{\M}$ is a valuation domain, we have either $a\vert b$ or $b\vert a$, in other words $\emptyset=V((a:b))\cap V((b:a))$.
\ebem
\pr Suppose that $(a:b)\notin \M$, then in $B_{\M}$ express $b=a.(a:b)^{-1}.(b:a)$. Now, suppose that in $B_{\M}$, $a\vert b$. Therefore $gcd^{B_{\M}}(a,b)=a$, so $(a:b)\in U_{\M}$, which exactly mean that $(a:b)\notin \M$.
\qed
\thm \label{ppel} Let $B$ be a B\'ezout domain. Then given any pp $\L_{B,V}$-formula $\phi(\bold z)$, there exists finitely many constructible subsets $C_{\phi,k}$ of $MSpec(B)$ such that for any $\M\in C_{\phi,k}$, $\phi(\bold z)$ is equivalent to a conjunction $\chi_{k}(\bold z)$ of atomic $\L_{B_{\M},V}$-formulas, in the classes $Mod_{B_{\M}}$ of all $B_{\M}$-modules, $k\in K$, $K$ finite.
\par In particular, any pp $\L_{B,V}$-formula $\phi(z)$ (in one variable) is equivalent (in $Mod_{B_{\M}}$) to a formula of the form $$z.a=0\;\&\;\bigwedge_{i=1}^n V_{\delta_i}(z.b_i),$$ for some $a, \;b_i\in B_{\M}$, $\delta_{i}\in \Gamma_{\M}^+$, $1\leq i\leq n$, with $\delta_{1}.v_{\M}(b_1)^{-1}>\cdots>\delta_{n}.v_{\M}(b_n)^{-1}$ and $v_{\M}(b_1)>\cdots >v_{\M}(b_n)$.
\ethm
\pr We will proceed by induction on the
number of existential quantifiers. We start with the existential $\L_{B,V}$-formula $\phi(\bold z):=\exists x_{n}\cdots\exists x_{1}\;\psi_{1}(x_{1},x_{2},\cdots,x_{n},\bold z)$, where $\psi_{1}(x_{1},x_{2},\cdots,x_{n},\bold z)$ is a
conjunction of atomic (c.a.) $\L_{B,V}$-formulas. 
We consider the innermost existential quantifier and the formula $\exists x_{1}\;\psi_{1}(x_{1},x_{2},\cdots,x_{n},\bold z)$. We want to find a finite covering of $MSpec(B)$ by constructible subsets $C_{1,j}$ and finitely many c.a. $\L_{B_{\M},V}$-formulas $\psi_{2,j}(x_{2},\cdots,x_{n},\bold z)$ such that for any $\M$ in $C_{1,j}$, $\exists x_{1}\;\psi_{1}(x_{1},x_{2},\cdots,x_{n},\bold z)$ is equivalent to $\psi_{2,j}(x_{2},\cdots,x_{n},\bold z)$ in $Mod_{B_{\M}}$.
\par We proceed inductively as follows. Set $C_{0,1}:=MSpec(B)$, $J_{0}:=\{1\}$. On the constructible subset $C_{\ell,j}$, $j\in J_{\ell}$, $n-1\geq \ell\geq 0$, we consider the existential formula $\exists x_{\ell+1}\;\psi_{\ell+1,j}(x_{\ell+1},\cdots,x_{n},\bold z)$, where $\psi_{\ell+1,j}$ is a c. a. $\L_{B_{\M},V}$-formula, and we show that 
there is a finite covering of $C_{\ell,j}$ by constructible subsets $C_{\ell+1,j'}$, $j'\in J_{\ell+1}$ such that for any $\M\in C_{\ell+1,j'}$, this formula $\exists x_{\ell+1}\;\psi_{\ell+1,j}(x_{\ell+1},\cdots,x_{n},\bold z)$ is equivalent in $C_{B_{\M}}$, to a formula $\psi_{\ell+2,j'}(x_{\ell+2},\cdots,x_{n},\bold z)$, where $\psi_{\ell+2,j'}$ is a c.a. $\L_{B_{\M},V}$-formula. 
\par At the last step, we obtain $J_{n}$, constructible subsets $C_{n,j}$, $j\in J_{n}$ and corresponding c.a. $\L_{B_{\M},V}$-formula $\psi_{n+1,j}(\bold z)$.
\par For ease of notation, we set $x:=x_{\ell+1}$ and $\bold y:=(x_{\ell+1},\cdots,x_{n},\bold z)$, $C_{\ell}:=C_{\ell,j}$ and $\psi_{\ell+1}(x,\bold{y}):=\psi_{\ell+1,j}(x,\bold{y})$.
\medskip
\par Among atomic $\L_{V}$-formulas, we have formulas of the form $V_{\delta}(x.r-u)$, where $r\in B$ and $\delta\in \Gamma^+$ that we will abbreviate as $x.r\equiv_\delta u$ (\lq\lq congruences relations\lq\lq). The outline of the proof is similar to \cite[Proposition 4.1]{BP}. 
\par First let us show that we can always
assume that we have at most one equation involving
$x.$ 
\par Indeed, consider $x.r_{0}=t_{0}(\bold y)\;\&\;x.r_{1}=t_{1}(\bold y)$.
For every $\M\in C_{\ell}$, either $v_{\M}(r_{0})\geq v_{\M}(r_{1})$, or $v_{\M}(r_{0})\leq v_{\M}(r_{1})$. W.l.o.g. assume we are in the second case. So, $r_{1}.r_{0}^{-1}\in B_{\M}$ and the above conjunction is equivalent to:
$x.r_{0}=t_{0}(\bold y)\;\&\;t_{0}(\bold y).r_{1}.r_{0}^{-1}=t_{1}(\bold y)$. So we will subdivise $C_{\ell}$ into two subsets according to whether  $v_{\M}(r_{0})\geq v_{\M}(r_{1})$ holds.
\par So we may reduce ourselves to consider c.a. $\L_{B_{\M},V}$-formulas $\psi_{\ell+1}(x,\bold{y})$ of the form $(\star)$, either 
$$x.r_{0}=t_{0}(\bold{y})\;\&\; \bigwedge_{i=1}^n\;V_{\delta_{i}}(x.r_{i}-
t_{i}(\bold{y}))\;\&\;
\theta(\bold{y}),\;\;{\rm or}$$
$$\bigwedge_{i=1}^n\;V_{\delta_{i}}(x.r_{i}-
t_{i}(\bold{y}))\;\&\;
\theta(\bold{y}),$$
where $r_{i}\in B_{\M}$, $\theta(\bold{y})$ is a c.a. $\L_{B_{\M},V}$-formula,
the $t_{i}(\bold{y})$, $0\leq i\leq n$, are $\L_{B_{\M}}$-terms, and
$\delta_{i}\in \Gamma_{\M}^+$. 
\par Consider $\exists
x\,\psi_{\ell+1}(x,\bold{y})$. It suffices to show that we can find finitely many c.a. $\L_{B_{\M},V}$-formulas and a finite covering of $C_{\ell}$ by finitely constructible subsets $C_{\ell+1,j}$ 
such that any such existential formula is equivalent to one of these formulas on $C_{\ell+1,j}$.

\par Before eliminating the existential quantifier, we make a series of reductions which lead us to break up $C_{\ell}$ into finitely many subsets according to whether certain valuational inequalities hold among the indices of the predicates or the coefficients occurring in the terms. To avoid too many indices after each reduction we rename $C_{\ell}$ each of the subsets we obtained from $C_{\ell}$.
\par First we examine whether
$v_{\M}(r_0)> v_{\M}(r_i)$, for $1\leq i\leq n$. Indeed, suppose 
that
$v_{\M}(r_{0})\leq v_{\M}(r_{i}),$ for some $i$, say $i=1$.
Then, we write $r_{1}=r_{0}.(r_{1}.r_{0}^{-1})\in B_{\M}$.
So, we replace in $(\star)$ the congruence relation $V_{\delta_1}(x.r_1-t_{1}(\bold{y}))$ by  
$V_{\delta_{1}}(t_{0}(\bold{y}).r_{1}.r_{0}^{-1}-
t_{1}(\bold{y}))$. So we are left with congruence relations $V_{\delta_i}(x.r_i-t_{i}(\bold{y}))$ with $v_{\M}(r_0)> v_{\M}(r_i)$.
\par We also break up $C_{\ell}$ according to finitely many subsets according to whether $v_{\M}(r_{i})<  \delta_{i\M}$, $1\leq i\leq n$.
Indeed suppose for instance that $v_{\M}(r_{i})\geq \delta_{i\M}$. Then we replace the congruence condition $x.r_{i}\equiv_{\delta_i}  t_{i}$ by $t_{i}\equiv_{\delta_i} 0$.
 \par Now let us order the set $\{\delta_{i\M}\cdot v_{\M}(r_{i}^{-1});\;1\leq i\leq n\}$. By re-indexing we may assume that $\delta_{1\M}\cdot v_{\M}(r_{1})^{-1}\geq \delta_{2\M}\cdot v_{\M}(r_{2})^{-1}\geq \cdots\geq \delta_{n\M}\cdot v_{\M}(r_{n})^{-1}.$
Again we break up $C_{\ell}$ into finitely many subsets according to which such conjunctions of inequalities hold.
\par First, we will assume that there is one
equation present in $\phi(x,\bold y)$ and we proceed by induction on the number of congruence conditions. If there are none, we replace the pp formula $x.r_{0}=t_{0}$ by $V_{v_{\M}(r_{0})}(t_{0})$. 
\par Consider the system ($1$):
$$(1) \  : \   x.r_0=t_0, \,
x.r_{1}\equiv_{\delta_1} t_{1}, \;
\cdots
,\, x.r_{n}\equiv_{\delta_n} t_{n} 
$$ 
\noindent with $v_{\M}(r_{0})> v_{\M}(r_{i}),$ 
 $t_0=t_{0}(\bold y), t_i=t_i(\bold y)$,  $1\leq i\leq n$.
\par We claim that in any $B_{\M}$-module $M$, with $\M\in C_{\ell}$, system ($1$) is equivalent to the
following system ($2$) :

$$(2) \ : \ x.r_{1}= t_{1},\;
t_{1}.r_{0}.r_{1}^{-1}\equiv_{\delta_{1}\cdot v(r_{0}.r_{1}^{-1})} t_{0}
,\;
x.r_{2}\equiv_{\delta_2} t_{2}, \;
\cdots 
, \; x.r_{n}\equiv_{\delta_n} t_{n}$$

$(1)\rightarrow (2)$\\
Let $x\in M$ satisfy (1). So $x.r_1=t_1+n$ for some $n\equiv_{\delta_1} 0$. Let $s_1\in B_{\M}$ be such that $v_{\M}(s_1)=\delta_1$. By definition of predicate $V_{\delta_{1}}$ and the assumption that $V_{\delta_{1}}(n)$,  there exists $u'\in M$ such that $u'.s_1=n$ and since $v_{\M}(r_1)<\delta_1$, we have $u'.(s_1.r_1^{-1}).r_1=n$ with $v_{\M}(s_1.r_1^{-1})>1$ (and so $s_1.r_1^{-1}\in B_{\M}$). Set $u:=u'.(s_1.r_1^{-1})$ and 
let $y:=x-u$. So we get that $y.r_1=t_1$ and $y.r_0=t_0-u.r_1.(r_0.r_1^{-1})=t_0-n.(r_0.r_1^{-1})=y.r_1.(r_0.r_1^{-1})=t_1.(r_0.r_1^{-1})$.
Therefore, $V_{\delta_{1}\cdot v(r_{0}.r_{1}^{-1})}(t_{0}-t_{1}.r_{0}.r_{1}^{-1})$.
Consider the other congruence conditions: if we replace $x$  by $y$, then for $i\geq 2$, $y.r_i=x.r_i-u'.(s_1.r_1^{-1}).r_i$ with $V_{\delta_i}(u'.(s_1.r_1^{-1}).r_i)$.\\
$(2)\rightarrow (1)$\\
Let $y$ satisfy $(2)$, namely $y.r_1=t_1$. Consider $y.r_0=y.r_1.(r_0.r_1^{-1})=t_1.(r_0.r_1^{-1})$.
Since $V_{\delta_{1}\cdot v(r_{0}.r_{1}^{-1})}(t_{0}-t_{1}.r_{0}.r_{1}^{-1})$ holds, we have $y.r_0\equiv_{\delta_{1}\cdot v(r_{0}.r_{1}^{-1})} t_0$.
Let $n\in M$ be such that $y.r_0=t_0+n$ with $n\equiv_{\delta_{1}\cdot v(r_{0}.r_{1}^{-1})} 0$. 
Let $s_1\in B_{\M}$ with $v_{\M}(s_1)=\delta_1$.  Since $V_{\delta_{1}\cdot v(r_{0}.r_{1}^{-1})}(n)$, there exists an element $n'\in M$ with $n'.s_1.(r_{0}.r_{1}^{-1})=n$, namely
$(n'.s_1.r_1^{-1}).r_0=n.$ Set $x=y-n'.s_1.r_1^{-1}$. Since $v_{\M}(s_1)=\delta_1$, $V_{\delta_{1}}(n'.s_1)$ holds, and similarly since we have $v_{\M}(s_1.r_1^{-1}.r_i)\geq \delta_{i}$, 
$V_{\delta_{i}}(n'.s_1.r_1^{-1}.r_i)$ holds, $i\geq 2$.
\medskip
\par Second, we will consider the case where there are only congruence relations in the system. We will either reduce to the previous case, making $x$ occurring in a non-trivial equation (see $(4)$). Let $M$ be any $B_{\M}$-module, with $\M\in C_{\ell}$.
\par Consider the following system ($3$):

$$(3) \  : \  
x.r_{1}\equiv_{\delta_1} t_{1}, \;
\cdots
,\, x.r_{n}\equiv_{\delta_n} t_{n} 
$$
 \par Since $v_{\M}(r_{1})< \delta_{1}$, we may replace the congruence condition by a divisibility condition. Indeed, $x.r_{1}=t_{1}+u$ with $u\equiv_{\delta_{1}} 0$. Let $s_{1}\in B_{\M}$ be such that $v_{\M}(s_{1})=\delta_{1}$ and let $u'\in M$ such that $u'.s_{1}=u$. 
 Set $y:=x-u'.s_{1}.r_{1}^{-1}$ and we check that for any $i\geq 2$, $V_{\delta_{i}}(u'.s_{1}.r_{1}^{-1}.r_{i})$.
 \medskip
\par So system ($3$) is equivalent to a
system ($4$) of the form:

$$(4) \ : \ x.r_{1}= t_{1},\;
x.r_{2}\equiv_{\delta_{2}} 
t_{2}
,\;
x.r_{3}\equiv_{\delta_3} t_{3}, \;
\cdots 
, \; x.r_{n}\equiv_{\delta_n} t_{n}.$$
Again we may assume that $v_{\M}(r_1)>v_{\M}(r_i)$, for all $i\geq 2$ (otherwise we may eliminate $x$ in the corresponding congruence relation). So we proceed as in the first case with one less congruence relation and we conclude by induction.
\medskip
\par Finally we consider a c.a. $\L_{B_{\M},V}$-formula with one free variable of the form $\bigwedge_{j\in J} x.a_j=0\;\&\;\bigwedge_{i\in I} V_{\delta_{j}}(x.b_{i})$. 
We proceed as in the beginning of the proof to reduce ourselves to at most one annihilator condition comparing the values $v_{\M}(a_i)$ in $\Gamma_{\M}^+$.
\par For the congruence conditions, as before we may assume that $v_{\M}(b_{i})< \delta_{i}$, otherwise we remove that the corresponding congruence relation. 
Then we compare in $(\Gamma(B_{\M}),\leq)$ the elements $\delta_{i}.v_{\M}(b_{i})^{-1}$ and so we assume that $\delta_{1}.v_{\M}(b_1)^{-1}\geq \cdots\geq \delta_{n}.v_{\M}(b_n)^{-1}$. 
Let $i\neq1$. First note that if $v_{\M}(b_1)\leq v_{\M}(b_i)$, we express $x.b_i=x.b_1.(b_1^{-1}.b_i)$. So if $V_{\delta_{1}}(x.b_1)$ holds, then $V_{\delta_{i}}(x.b_i)$ holds.
Now assume that $\delta_{1}.v_{\M}(b_1)^{-1}=\delta_{i}.v_{\M}(b_i)^{-1}$, $v_{\M}(b_1)>v_{\M}(b_i)$ and $V_{\delta_i}(x.b_i)$. Then since $x.b_1=x.b_i.b_1.b_i^{-1}$, we have $V_{\delta_{1}}(x.b_{1})$.
So, proceeding in a similar way for all indices, we may assume that the congruence conditions are of the form: $\bigwedge_{i=1}^n V_{\delta_{j}}(x.b_{i})$ with $\delta_{1}.v_{\M}(b_1)^{-1}> \cdots> \delta_{n}.v_{\M}(b_n)^{-1}$ and $v_{\M}(b_1)>\cdots>v_{\M}(b_n)$. \qed
\bem \label{cons} From the proof of the above theorem, we see that the constructible subsets $C_{\phi,k}$ occurring in the statement are of the form $V(a)^c\cap V(b)$, for $a, b\in B$ that can be obtained from the formula $\phi$ and the operations $\cdot$ and $gcd$ (in $B$).
\ebem

\section{Axiomatization}
\par Since we are interested in decidability results for theories of modules, we will axiomatize the theories of modules that we will consider. We will start with considering modules over any B\'ezout ring and then we will apply the results of the preceding section on valuation domains.
\bigskip 
\par For a ring $R$, when one considers the decidability problem for the theory of $R$-modules, it is reasonable to assume from the start, certain effectivity properties of the ring $R$. In particular one can ask: assuming that the theory of $R$-modules is decidable, which effectivity properties does it imply on the ring operations? One usually assume that $R$ is {\it effectively given}. This notion has been discussed in length in  \cite{PPr}, \cite[Chapter 17]{Pr} and specifically for
valuation domains in: \cite[Section 3]{PPT}, \cite[Definition 3.1]{Gr}). In particular if a valuation domain $A$ is effectively given, there is an algorithm which given $a,\,b\in A$ decides whether $a\vert b$ \cite[Remark 3.2]{PPT}.
\medskip
\defn\label{EF} Let $B$ be a  B\'ezout domain, we will call assumption $(EF)$ on $B$ the following: $B$ is a countable ring, it can be enumerated as $(r_{n}:\;n\in \omega)$ in such a way, there are algorithms to perform the following operations: given $a,\;b\in B$, produce $a+b$, $-a$, $a.b$, decide whether $a=b$ or not and the relation $\{(a,b)\in R^{2}:\;a\vert b\}$ is recursive. 
\edefn
\medskip
\par This implies (as in \cite[Section 3]{PPT}) there is an algorithm which decides whether an element $a\in B$, decide whether $a$ is invertible (i.e. $a\in U$) and if yes produce $a^{-1}$. Therefore given any coset $a.U$, there is an algorithm which chooses a representative (for instance the element with the smallest index in the given enumeration). There are also algorithm which given $a,\;b\in B^{\star}$, produces:
\par  $gcd(a,b)$ (the algorithm enumerates the elements $a.r_{n}+b.r_{m}$ and checks whether it divides $a$ and $b$),
\par $(a:b)$ (the algorithm checks whether $a=gcd(a,b).r_n.u$ for some $u\in U$).
\subsection{Abelian structures revisited}
\par Recall that the language $\L_{V}$ has been defined as the expansion of the language $\L$ of $B$-modules together with a set of unary predicates that we index by the submonoid of the positive elements of the group of divisibility $\Gamma$ of $B$, namely $\L_{V}:=\L\cup\{V_{\gamma}; \gamma\in\Gamma^+\}$. 
\defn \label{abelian}
Let $T_{B,V}$ be the $\L_{V}$-theory, consisting of the $\L$-theory of $B$-modules together with:
\begin{enumerate}
\item[$(1)_{V}$] $V_{\delta}(0)$, for each $\delta\in \Gamma^+$,
\item[$(2)_{V}$] $\forall m_{1}\;\forall m_{2}\;(V_{\delta_{1}}(m_{1})\,\&\,V_{\delta_{2}}(m_{2}))\rightarrow V_{\delta_{1}\wedge\delta_{2}}(m_{1}+m_{2})$, for every $\delta_{1},\;\delta_{2}\in \Gamma^+$,
\item[$(3)_{V}$] $\forall m\;(V_{\delta}(m)\rightarrow V_{\delta\wedge \mu}(m))$, for any $\delta,\;\mu\in \Gamma^+$,
\item[$(4)_{V}$] $\forall m\;V_{\delta}(m)\rightarrow V_{\delta\cdot v(a)}(m.a)$, for each $a\in B^{\star}$, $\delta\in \Gamma^{+}$.
\end{enumerate}
\edefn
 \par When the context is clear, we will simply use the notation $T_{V}$ (instead of $T_{B,V}$). Given a $B$-module $M$, we denote by $M_V$ its expansion as an $\L_V$-structure, namely $M$ together with a family a submodules $V_\delta(M)$, $\delta\in \Gamma^+$.
\bem\label{property3} From the above axioms, we deduce easily the following properties, letting $M_V\models T_{B,V}$:
$\;$
\begin{enumerate}
\item each $V_{\delta}(M)$ is a subgroup and it is a $B$-submodule of $M$,
\item if $u\in B$ is an invertible element e.g. $u\in U$, then for any $\delta\in \Gamma$, $V_{\delta}(m.u)\leftrightarrow V_{\delta}(m)$ (this is due to the fact that $v(u)=1$),
\item suppose $v(a)=v(b)$ with $a,\;b\in B^{\star}$ and that $V_{\delta.v(b)}(m.b)$, then $V_{\delta.v(b)}(m.a)$.
\end{enumerate}
\ebem
\medskip
\par We will consider a subclass of the class of abelian structures we just defined, namely those which satisfy in addition:
$$(5)_{V}\;\;\;\;\;\;\;\;\;\;\; \forall m\in M\;V_{1}(m),$$
 together with the following divisibility scheme: 
$$(6)_{V,div}\;\;\;\;\;\;\;\;\;\;\;\forall m\;\exists n \;(V_{v(a)}(m)\;\rightarrow n.a=m), {\rm for\;each}\; a\in B^{\star}.$$
\defn\label{div} We will denote by $T_{B,V,div}$ the theory $T_{B,V}$ together with $(5)_V$ together with the divisibility axioms scheme  $(6)_{V,div}$.
\edefn
\bem $\;$
\begin{enumerate}
\item The theory $T_{B,V,div}$ is consistent.\\
Indeed, the ring $B$ itself can be expanded to a model of $T_{B,V,div}$. Define $V_{\gamma}(B)=\{b\in B: v(b)\geq \gamma\}$, $\gamma\in \Gamma^+$. If $\gamma:=v(a)$, $a\in B^{\star}$, $V_{\gamma}(B)=B.a$.
Then $(B,(Ba)_{a\in B^{\star}})$ is a model of $T_{B,V,div}$.  By definition of the $\ell$-valuation on $B$, we have that for $a,\;b\in B^{\star}$, $v(b)\geq v(a)$ iff $b.a^{-1}\in B$ and so $B$ will satisfy axiom $(5)_V$ together with axiom schemes $(4)_V$ and $(6)_{V,div}$. In fact, $B$ satisfies a stronger form of axiom scheme $(4)_{V}$, as we will see below. 
\item In fact, $T_{B,V,div}$ is what is called {\it an expansion by definitions} of the theory of $B$-modules. 
Explicitly, it means that given any $B$-module $M$, we can expand it to a model of $T_{B,V,div}$ and in any model of $T_{B,V,div}$, the new predicates $V_{\gamma}$ are definable in the language of $B$-modules. \\								
Indeed, given any $B$-module $M$, we set $V_\delta(M):=\{m\in M:\;\exists n\;m=n.a\}$ with $v(a)=\delta$, $a\in B^{\star}$.  
To check that this expansion satisfies all the axioms $(1)_V$ up to $(6)_{V,div}$ is rather straightforward. Let us check for instance $(2)_{V}$. Let $\delta_1=v(a_1)$ and $\delta_2=v(a_2)$ and write $a_1=a'.gcd(a_1,a_2)$, $a_2=a''.gcd(a_1,a_2)$. Assume $V_{\delta_{1}}(m_{1})\,\&\,V_{\delta_{2}}(m_{2})$ holds. So for some $n_1,\;n_2$ we have $m_1=n_1.a_1$ and $m_2=n_2.a_2$ and $m_1+m_2=(n_1.a'+n_2.a'').gcd(a_1,a_2)$. Therefore $V_{\delta_1\wedge\delta_2}(m_1+m_2)$.
\par By definition of the $\ell$-valuation on $B$, we have that for $a,\;b\in B^{\star}$, $v(b)\geq v(a)$ iff $b.a^{-1}\in B$. Using this, one can easily deduce axioms  $(3)_V$, $(5)_V$ and $(6)_{V,div}$.
\par Conversely, if $M_{V}\models T_{B,V,div}$, then the subgroup $V_{v(a)}(M)$, $a\in B^{\star}$ is definable by the $\L_{B}$-formula: $\exists x\;y=x.a$. It follows from $(5)_{V}$ and $(6)_{V,div}$. 
\item
Let $M_{V}\models T_{B,V,div}$ and assume $M$ is a torsion-free $B$-module. Let $m\in M$ and suppose that $V_{\delta.v(b)}(m.b)$ holds with $b\in B^{\star}$, $\delta\in \Gamma^{+}$. 
Then let us show that $V_{\delta}(m)$ holds. Let $a\in B^{\star}$ be such that $v(a)=\delta$.
 By axiom $(6)_{V,div}$, for some $n\in M$, $m.b=n.a.b$. So $(m-n.a).b=0$. Since $M$ is torsion-free, $m= n.a$ and so $V_{\delta}(m)$ holds. 
\end{enumerate}
\ebem
\subsection{The case of valuation domains}
\par As recalled in the introduction, there are general results on decidability of the theories of modules over valuation domains \cite{PPT}, \cite{Gr}. For instance, G. Puninski, V. Puninskaya and C. Toffalori proved that if $A$ is a valuation domain satisfying $(EF)$ with infinite residue field and archimedean densely ordered value group, then the theory of all $A$-modules is decidable \cite[Theorem 6.2]{PPT}. Then L. Gregory removed these two assumptions and proved that for an effectively given valuation domain $A$, if the prime radical relation is recursive, then the theory of $A$-modules is decidable \cite[Theorem 7.1]{G}. In addition she proved for any effectively given commutative ring $R$, if the theory of the $R$-modules is decidable then the prime radical relation is recursive. 

\par In order to be self-contained, we will present here a direct proof of their decidability result, under the hypothesis of infinite residue field. As said in the introduction, this proof is also more algebraic than the one presented in \cite{PPT}.
It will also be a key step in the decidability result for B\'ezout domains in the next section. The hypothesis on the residue field is satisfied in all good Rumely domains containing the prime field $\F_p$ (see section \ref{App}).
\medskip
\par In this subsection, $A$ will denote a valuation domain, $Q(A)$ its fraction field and $\Gamma:= \Gamma(A)$ its value group. Recall that a fractional ideal of $A$ is an additive subgroup of $Q(A)$ which is closed under multiplication by $A$.
We will use that a pure-injective indecomposable module over a valuation domain is the pure hull $\overline{I/J}$ of a module of the form $I/J$, where $I, J$ are two fractional ideals of $A$ and that if a type is realized in $\overline{I/J}$, it is already realized in $I/J$ \cite[section 5]{Z}.
\par Note that $v(J)$ is a subset of $\overline{\Gamma}$ which is upward closed (in particular if $x\in I\setminus J$, then $v(x)<v(J)$). We define the predicates $V_{\delta}$ on $I/J$ using axiom $(6)_{V,div}$. Namely  let $u\in I\setminus J$, then $V_{\delta}(u+J)$ holds if there exists $s\in A$ with $v(s)\geq \delta$ and $y\in I$ such that $u-y.s\in J$. 
We have that $v(u)\geq min\{y.s,u-ys\}$, since $u\notin J$, $y.s\notin J$ and so $v(ys)<v(J)$ and $v(u)=v(y.s)$, so $v(u)\geq\delta$ and $s\vert u$ (in $I$).
\nota Let $E$ be a subset of $\Gamma$ which is upward closed. Let $^*\Gamma$ be an elementary extension of $\Gamma$ which is $\vert \Gamma\vert^+$-saturated. We consider the partial type with parameters in $\Gamma$, $p(x):=\{\delta<x\leq \gamma:\;\delta\in \Gamma\setminus E,\;\gamma\in E\}$ and we denote by $\inf E$ a realisation of $p(x)$ in $^*\Gamma$.
\enota
\par Let $E$ be a subset of $\Gamma$ which is upward closed and assume that $E$ has no minimum in $\Gamma$. Then any non empty interval in $^*\Gamma$ of the form $]\inf E\;\gamma]$, where $\gamma\in \Gamma$ has infinite intersection with $\Gamma$.
\cor \label{dec-val} Assume that $A$ satisfies $(EF)$, and that the quotient $A/\M$ of $A$ by its maximal ideal $\M$, is infinite, then $T_{A}$ is decidable.
\ecor 
\pr 
 Since the theory $T_{A,V,div}$ is recursively enumerable, in order to prove its decidability we need to show that we can enumerate the set of sentences false in some element of $T_{A,V,div}$ (or taking the negation, true in some element of $T_{A,V,div}$). 
 \par By the Baur--Monk pp elimination theorem, it suffices to consider boolean combinations of invariant sentences: $(\phi/\psi)\geq n$, $n\in \N$, where $\phi(x), \psi(x)$ are pp $\L_{V}$-formulas (see subsection \ref{prelim2}). Note that since the quotient of $A$ by its maximal ideal $\M$ is infinite, then $(\phi/\psi)>1$ implies that  $(\phi/\psi)$ is infinite.
\par Since a disjunction of formulas is true whenever one of them is true, we may only consider conjunctions of formulas of the form $\si:=\bigwedge_{i\in I} (\phi_i/\psi_i)>1\;\&\;\bigwedge_{j\in J} (\chi_j/\xi_j)=1$. Moreover suppose we find for each $i\in I$, a model $M_i$ of $T_{A,V,div}$ satisfying $(\phi_i/\psi_i)>1\;\&\;\bigwedge_{j\in J} (\chi_j/\xi_j)=1$, then we can form the direct sum of the $M_i's$, $i\in I$ and get a model $M$ of $T_{A,V,div}$ satisfying $\si$.
So the sentences we need to consider are of the following form: $(\phi_i/\psi_i)>1\;\&\;\bigwedge_{j\in J} (\chi_j/\xi_j)=1$. Furthermore note that if we cannot find a model $M_{V}$ of $T_{A,V,div}$ where $\si$ holds, it exactly means that in the Ziegler spectrum the basic open set $[\phi_i/\psi_i]$ is included in the union: $\bigcup_{j\in J} [\chi_{j}/\xi_{j}]$.
 \par By Theorem \ref{ppel}, we may reduce ourselves to only consider pp $\L_{V}$-formulas of the form: $x.a=0\;\&\;\bigwedge_{i\in I} V_{\delta_i}(x.b_i)$ with $a,\;b_i\in A^*$, $\delta_{i}\in \Gamma(A)$, $I$ finite, and or equivalently, $x.a=0\;\&\;\bigwedge_{i\in I} ((x.b_i=0)+V_{\delta_i.v(b_i)^{-1}}(x))$, assuming that $v(b_i)<\delta_i$ (otherwise we may delete the corresponding congruence condition). In order to write that formula in the same way when $I$ is empty, we will express $x.a=0$ by $V_{+\infty}(x.a)$.
 \par Furthermore, note that the quantifier elimination procedure described in Theorem \ref{ppel} is effective. In the course of the proof, we had to decide
whether $v(a_{1})<v(a_{2})$, $a_{1},\;a_{2}\in A^{\star}$. This is equivalent to decide whether 
$a_{1}\vert a_{2}$ and $a_{2}\nmid a_{1}$. By hypothesis (EF), we can do that in an effective way in $A$. 
 \medskip
\par As in \cite{PPT}, we use the duality functor introduced by M. Prest on the lattice of pp formulas (see subsection \ref{prelim2}) in order to simplify the form of the pairs of pp formulas we need to consider. Let $\phi(x)$ be a pp formula, then $D(D(\phi))\leftrightarrow \phi$. Assuming that $D(\phi(x))\leftrightarrow \bigwedge_i(c_i\vert x+x.d_i=0)$, we get $\phi(x)\leftrightarrow \sum_i D(c_i\vert x+x.d_i=0)\leftrightarrow \sum_i\,(d_i\vert x\;\&\;x.c_i=0)$.
Finally we note that if $\phi\leftrightarrow \sum_i \phi_i$ and $\psi\leftrightarrow \bigwedge_j \psi_j$, then $[\phi/\psi]=[\bigcup_{i,j} (\phi_i/\psi_j]$. Since the formula $a\vert x+x.d=0$ is equivalently $V_{v(a.d)}(x.d)$, by making the same abuse of notation as above by allowing the possibility to have $V_{+\infty}(x.d)$, we get the following Claim.
\cl {\rm \cite[Section 5]{PPT}}\label{1} We only need to consider basic open sets in the Ziegler spectrum of the form $[b\vert x\;\&\;x.c=0/V_{v(a.d)}(x.d)]$.
\ecl
\cl {\rm \cite[Corollary 4.3]{PPT}} \label{2} $[x.b=0\;\&\;V_{\delta_{1}}(x)/V_{\delta_{2}}(x.c)]\neq \emptyset$ if and only if ($v(b)>v(c)$ and $\delta_1.v(c)<\delta_2$).
\ecl
\pr $(\leftarrow)$ Let $s\in A$ be such that $v(s)=\delta_1$. Consider the $A$-module $M:=A/A.bs$; the element $x:=s+A.bs$ belongs to $ann(b)\cap V_{\delta_1}(M)$. By the way of contradiction, assume that $V_{\delta_2}(s.c+A.bs)$. We have that $v(s.c)=\delta_1.v(c)<min\{\delta_1.v(b),\delta_2\}$, a contradiction.
\par $(\rightarrow)$ Let $N$ be an $A$-module such that there exists $m\in N$ such that $m.b=0,$ $V_{\delta_1}(m)$ and $\neg V_{\delta_2}(m.c)$.
If $v(b)\leq v(c)$, then $m.c=0$ and so $V_{\delta_2}(m.c)$. Now assume for a contradiction that $\delta_1.v(c)\geq \delta_2$. Since $V_{\delta_1}(m)$, we have that $V_{\delta_2}(m.c)$, contradicting the assumption on $m$. \qed
\bigskip
\par Now, let us consider the two open sets: $[x.b=0\;\&\;V_{\delta_{1}}(x)/V_{\delta_{2}}(x.c)]$,\\$[x.b'=0\;\&\;V_{\delta_{1}'}(x)/V_{\delta_{2}'}(x.c')]$, with $b,\,c,\;b',\;c'\in A^{*}$, $\delta_{1},\;\delta_{2}, \delta_{1}',\;\delta_{2}'\in \Gamma^{+}$ and\\ $v(b)>v(c)$, $v(b')>v(c')$, $\delta_1.v(c)<\delta_2$, $\delta_1'.v(c')<\delta_2'$. 
\cl {\rm (See also \cite[Proposition 4.5]{PPT})} \label{3} Under the above assumptions, we have \\ 
$[x.b=0\;\&\;V_{\delta_{1}} (x)/V_{\delta_{2}}(x.c)]\subset [x.b'=0\;\&\;V_{\delta_{1}'}(x)/V_{\delta_{2}'}(x.c')]$ holds if and only if \\
$\delta_2.v(c)^{-1}\leq \delta_2'.v(c')^{-1}$ and
$\delta_2.v(c)^{-1}.v(b)\leq \delta_2'.v(c')^{-1}.v(b').$
\ecl
\pr $(\rightarrow)$
\par $(i)$ By the way of contradiction, assume that $\delta_{2}'.v(c')^{-1}<\delta_2.v(c)^{-1}$. Choose $I,\;J$ two fractional ideals of $Q(A)$, $J\subset I$ such that $\min(v(I))=\delta_2'.v(c')^{-1}$ and $\min(v(J))=\delta_2.v(c)^{-1}$. 

Let us show that $\overline{I/J}$ belongs to the first open set but not to the second one.
Let $u\in I$ with $v(u)\geq \delta_{1}$ and $v(u).v(b)\in v(J)$; by choice of $I,\;J$ we have that $v(u)\geq \delta_1$ and for any $m\in J$, $v(u.c+m)<\delta_2$. However all elements of $I$ have valuation $\geq \delta_2'.v(c')^{-1}$.
\par $(ii)$ Now, let us show that if $\delta_2.v(c)^{-1}.v(b)>\delta_2'.v(c')^{-1}.v(b')$, then we get a contradiction.
Choose a fractional ideal $J$ such that $\max\{\delta_2'.v(c')^{-1}.v(b'),\delta_2\} \leq \min(v(J))<\delta_2.v(c)^{-1}.v(b)$. This is feasible since $\delta_2< \delta_2.v(c)^{-1}.v(b)$.
Let $u\in A$ be such that $\delta_1\leq v(u)<\delta_2.v(c)^{-1}$ and $v(u.b)\in v(J)$, equivalently $\min(v(J)).v(b)^{-1}\leq v(u)$. Moreover, since $\delta_2\leq v(J)$,  we have that $v(u.c+m)<\delta_2$, for any $m\in J$.
However, any element $\tilde u\in I\setminus J$ such that $v(\tilde u.c')<\delta_2'$ has the property that $\tilde u.b'\notin J$, a contradiction.
\medskip
\par $(\leftarrow)$ Now take any pure-injective indecomposable module belonging to the first pair. As already recalled, this module is of  the form $\overline{I/J}$, where $I,\;J$ are two fractional ideals. Since $\overline{I/J}$ belongs to the first pair, then there exists $u\in I\setminus J$ with $v(u)\in [\delta_1\;\delta_2.v(c)^{-1}[$, $\neg V_{\delta_2}(u.c+J)$ and such that $v(u.b)\in v(J)$. Note that this implies that $u.c\notin J$ and since $v(J)$ is upward closed, it implies that $v(u.c)<v(J)$, which implies that $\delta_1.v(c)\leq v(u.c)<v(J)$. 
\par Now we look for an element $u'\in I$, such that $v(u')\geq \delta_1'$ and it belongs to the interval $[\inf v(J)v(b')^{-1}\;\delta_2'.v(c')^{-1}[$. By hypothesis the interval $[\inf v(J)v(b)^{-1}\;\delta_2.v(c)^{-1}[$ is non trivial. This is equivalent to  $[\inf v(J)\;\delta_2.v(c)^{-1}v(b)[\neq \emptyset$. Since $\delta_2'.v(c')^{-1}.v(b')\geq \delta_2.v(c)^{-1}.v(b)$, $[\inf v(J)v(b')^{-1}\;\delta_2'.v(c')^{-1}[\neq \emptyset$. The interval $[\delta_1'\;\delta_2'.v(c')^{-1}[$ is non trivial, as well as $[\inf v(I)\;\delta_2'.v(c')^{-1}[$ (since $[\inf v(I)\;\delta_2.v(c)^{-1}[\neq \emptyset$ and $\delta_{2}.v(c)^{-1}\leq \delta_2'.v(c')^{-1}$). So the intersection of these three intervals is non empty.
It remains to check that $\neg V_{\delta_2'}(u.c'+J)$, namely for all $m\in J$, $v(u'.c'+m)<\delta_2'$. 
\par We always have that $\inf v(J).v(b')^{-1}<\inf v(J). v(c')^{-1}$. Either, $\inf v(J). v(c')^{-1}<\delta_2'.v(c')^{-1}$, in which case we replace the interval  $[\inf v(J)v(b')^{-1}\;\delta_2'.v(c')^{-1}[$ by \\$[\inf v(J)v(c')^{-1}\;\delta_2'.v(c')^{-1}[$
or $\inf v(J)\geq \delta_2'$. In that last case, any element $m\in J$ will have the property that $v(u'.c'+m)<\delta_2'$ provided that $v(u'.c')<\delta_2'$.
\qed
\smallskip
\par Before  considering the general case, let us consider the case when an open set in the Ziegler spectrum is included in the union of two open subsets.
We assume that each of the open sets is non trivial (see Claim \ref{2}).
\cl \label{4} $[x.b=0\;\&\;V_{\delta_{1}}(x)/V_{\delta_{2}}(x.c)]\subset$\\
$[x.b'=0\;\&\;V_{\delta_{1'}}(x)/V_{\delta_{2'}}(x.c')]\cup  [x.b''=0\;\&\;V_{\delta_{1''}}(x)/V_{\delta_{2''}}(x.c'')]$ if and only if,\\ either 
($\delta_2.v(c)^{-1}\leq \delta_2'.v(c')^{-1}$ and $\delta_2.v(c)^{-1}.v(b)\leq \delta_2'.v(c')^{-1}.v(b')$), or\\
($\delta_2.v(c)^{-1}\leq \delta_2''.v(c'')^{-1}$ and $\delta_2.v(c)^{-1}.v(b)\leq \delta_2''.v(c'')^{-1}.v(b'')$).
\ecl
\pr $(\rightarrow)$ Suppose otherwise. By symmetry, it suffices to consider the following cases:
\par (i) $\delta_2''.v(c'')^{-1}<\delta_2.v(c)^{-1}$ and $\delta_2'.v(c')^{-1}.v(b')<\delta_2.v(c)^{-1}.v(b)$.
\par Choose a fractional ideal $J$ such that $\max\{\delta_2'.v(c')^{-1}.v(b'),\delta_2\} \leq \min(v(J))<\delta_2.v(c)^{-1}.v(b)$
and a fractional ideal $I$ such that $\min v(I)=\max\{\delta_{1},\delta_{2}''.v(c'')^{-1}\}$. First let us check that $\overline{I/J}$ belongs to $[x.b=0\;\&\;V_{\delta_{1}}(x)/V_{\delta_{2}}(x.c)]$. Since $\min(v(J))<\delta_2.v(c)^{-1}.v(b)$, there is an element $u\in I$ such that $v(u.c)<\delta_{2}$ but $u.b\in J$ (and by choice of $I$, $v(u)\geq \delta_{1}$).
\par But any  element $u\in I$ will have the property that $v(u.c'')\geq \delta_{2}''$ and any element $u'\in I\setminus J$ with $v(u'.c')<\delta_{2}'$ will have the property that $u'.b'\notin J$.
\par(ii) $\delta_2.v(c)^{-1}>\max\{ \delta_2'.v(c')^{-1}, \delta_2''.v(c'')^{-1}\}$. Choose a fractional ideal $I$ such that $\min v(I)=\max\{\delta_{1},\delta_2'.v(c')^{-1}, \delta_2''.v(c'')^{-1}\}$ and a fractional ideal $J$ with $\delta_{2}\leq \min(J)<\delta_{2}.v(c)^{-1}.v(b)$. We similarly check that $\overline{I/J}$ belongs to $[x.b=0\;\&\;V_{\delta_{1}}(x)/V_{\delta_{2}}(x.c)]$.
But, any $u\in I\setminus J$ will have the property that $V_{\delta_{2}'}(u.c')$ and $V_{\delta_{2}''}(u.c'')$, but there is $u\in I$ such that $v(u)<\delta_{2}.v(c)^{-1}$ and $u.b\in J$.
\par (iii) $\delta_2.v(c)^{-1}.v(b)> \max\{ \delta_2'.v(c')^{-1}.v(b'), \delta_2''.v(c'')^{-1}.v(b'')\}$. Choose a fractional ideal $J$ such that $\max\{ \delta_2'.v(c')^{-1}.v(b'), \delta_2''.v(c'')^{-1}.v(b'')\} \leq \min v(J)<\delta_2.v(c)^{-1}.v(b)$. Choose a fractional ideal $I$ with the property that $\min(I)=\delta_{1}$. Again, it is easily checked that $\overline{I/J}$ belongs to $[x.b=0\;\&\;V_{\delta_{1}}(x)/V_{\delta_{2}}(x.c)]$.
\par But no $u\in I\setminus J$ such that $\neg V_{\delta_{2}'}(u.c')$ has the property that $u.b'\in J$. Similarly $\neg V_{\delta_{2}''}(u.c'')$ implies that $u.b''\notin J$.

\par $(\leftarrow)$ This direction is clear using the previous claim.
\qed
\medskip
\par Finally we will show that to decide whether 
a basic open set in the Ziegler spectrum is included in a finite union of basic open subsets reduces to divisibility conditions on the elements of the ring, which we can decide by assumption (EF).
As before, we assume that each of the basic open sets is non trivial (see Claim \ref{2}).
\cl \label{5} $[x.b=0\;\&\;V_{\delta_{1}}(x)/V_{\delta_{2}}(x.c)]\subset \bigcup_{\ell\in L} [x.b_\ell=0\;\&\;V_{\delta_{1\ell}}(x)/V_{\delta_{2\ell}}(x.c_{\ell})]$ if and only if 
$\bigvee_{\ell\in L}\;(\delta_2.v(c)^{-1}\leq \delta_{2\ell}.v(c_{\ell})^{-1}\;\;{\rm and}\;\;\delta_2.v(c)^{-1}.v(b)\leq \delta_{2\ell}.v(c_{\ell})^{-1}.v(b_{\ell})).$
\ecl
\pr $(\rightarrow)$ Suppose otherwise, namely that either 
\par (i) we can partition $L$ into two non-empty subsets $L', L''$ with
$\max_{\ell\in L''} \delta_{2\ell}.v(c_{\ell})^{-1}<\delta_2.v(c)^{-1}\leq \min_{\ell\in L'} \delta_{2\ell}.v(c_{\ell})^{-1}$ and 
$\max_{\ell\in L'} \delta_{2\ell}.v(c_{\ell})^{-1}.v(b_{\ell})<\delta_2.v(c)^{-1}.v(b)$.
In this case, we choose a fractional ideal $J$ such that $\max\{\max_{\ell\in L'} \delta_{2\ell}.v(c_{\ell})^{-1}.v(b_{\ell}),\delta_2\} \leq \min(v(J))<\delta_2.v(c)^{-1}.v(b)$
and a fractional ideal $I$ such that \\$\min v(I)=\max\{\delta_{1},\max_{\ell\in L''} \delta_{2\ell}.v(c_{\ell})^{-1}\}$. First let us check that $\overline{I/J}$ belongs to $[x.b=0\;\&\;V_{\delta_{1}}(x)/V_{\delta_{2}}(x.c)]$. Since $\min(v(J))<\delta_2.v(c)^{-1}.v(b)$, there is an element $u\in I$ such that$v(u.c)<\delta_{2}$ but $u.b\in J$ (and by choice of $I$, $v(u)\geq \delta_{1}$).
\par But any  element $u\in I$ will have the property that $v(u.c_\ell)\geq \delta_{2\ell}$, for $\ell\in L''$ and any element $u'\in I\setminus J$ with $v(u'.c_\ell)<\delta_{2\ell}$ will have the property that $u'.b_\ell\notin J$, $\ell\in L'$.
\par (ii) $\max_{\ell\in L} \delta_{2\ell}.v(c_{\ell})^{-1}<\delta_2.v(c)^{-1}$, then we proceed as in Claim \ref{4} (ii).
\par (iii) $\max_{\ell\in L} \delta_{2\ell}.v(c_{\ell})^{-1}.v(b_{\ell})<\delta_2.v(c)^{-1}.v(b)$, then we proceed as in Claim \ref{4} (iii).
\par $(\leftarrow)$ This direction is clear using Claim \ref{3}.
\par This ends the proof of the Claim and the proof of the Corollary. 
\qed
 \qed
\medskip
\bem \label{Zg} Let $T$ be a theory of $R$-modules, where the invariant sentences are of the form $(\phi/\psi)>1$. The discussion above showed that an equivalent formulation of whether a sentence holds in some $R$-module is asking 
whether a basic open set is included in a given finite union of other basic open sets in the closed subset of Ziegler spectrum of $R$, consisting of models of $T$ (see for instance \cite[section 6]{PPT}).  
\par When $B$ is a B\'ezout domain, each point of the Ziegler spectrum is an indecomposable pure-injective $B$-module (and so a $B_{\M}$-module, where $B_{\M}$ denotes the localization of $B$ at $\M$, for some $\M\in MSpec(B)$)
and a basic open set is the set of points in the Ziegler spectrum where the index the two pp definable subgroups is nontrivial. 
By the discussion above, in case $B_{\M}$ is a model of $T_{B_{\M},V,div}$ we reduce ourselves to consider open sets in the Ziegler spectrum of the form $[x.b=0\;\&\;V_{\delta_{1}}(x)/V_{\delta_{2}}(x.c)]$.
\ebem

\section{Feferman-Vaught theorem for B\'ezout domains}\label{FV}
\par Let $B$ be a B\'ezout domain and let $\Gamma(B)$ be its group of divisibility. 
\par S. Garavaglia \cite{G} showed that any $B$-module $M$, can be embedded in a direct sum of modules over the localizations $B_{\M}$, $\M$ varying in the space $MSpec(B)$ of maximal ideals of $B$, and this embedding is elementary (i.e. respects pp formulas). Even though we could have directly applied his result, we will present here a slight generalization for abelian structures (see Proposition \ref{gar} below)). 

\medskip
\subsection{Localizations}
\par Let us review basic definitions on localizations by maximal ideals of both the ring and the module \cite[Chapter 9]{GW}.
\par Let $\M\in MSpec(B)$. Let $M$ be a $B$-module and let $M_{\M}$ be the localization of $M$ by $\M$. There is an embedding of $M$ into the direct product $\prod_{\M\in MSpec(B)} M_{\M}$ (as a $B$-module) and S. Garavaglia showed that this embedding is elementary \cite[Theorem 3]{G} (namely respects pp formulas).
We want to extend this result when $M$ is viewed as an $\L_{V}$-structure, namely not only as a $B$-module but endowed with a distinguished lattice of submodules; to that end, we will use the following description of the localizations $M_\M$ and of the embedding of $M$ into the direct product $\prod_{\M\in MSpec(B)} M_{\M}$.
\par Recall that $M_{\M}$ is also a $B_{\M}$-module and that one can view $M_{\M}$ as the module of fractions $M\otimes B_{\M}$ of $M$ with respect to the multiplicative set: $B\setminus \M$ (\cite[Proposition 9.14]{GW}).   
\par Let $E\subset B$, denote by $ann_{E} M:=\{m\in M:\;\exists r\in E\;m.r=0\}$. 
For $m\in M^{\star}$, let $Ann(m):=\{r\in B:\;m.r=0\}$; it is a proper ideal of $B$ and so it is included in a maximal ideal $\M_{0}$ of $B$. Therefore $m\notin ann_{B\setminus \M_{0}} M$. 
\par We can embed $M$ into $\prod_{\M} M_{\M}$ as follows. Set $m_{\M}:=m+ann_{B\setminus \M} M$, with $m\in M$. The map sending $m$ to $(m_{\M})_{\M\in MSpec(B)}$ is injective by the above and clearly a morphism of $B$-modules. 
\par Now we consider the expansion of $M$ to the abelian structure $M_{V}$ as defined in Definition \ref{abelian}. It induces the following abelian structure on $(M_{\M})_{V}$ by setting 
$V_\delta(m_\M)$, whenever there exists $n\in m+ann_{B\setminus \M} M$ such that $V_\delta(n)$, where $\delta\in \Gamma$. 
\lmm Let $M$ be a $B$-module and let $M_V$ be its expansion as an abelian $\L_V$-structure. Let $\M\in MSpec(B)$.
Then there is a morphism of $\L_V$-structures sending $m\in M_V$ to $m_{\M}\in (M_{\M})_{V}$.
\elmm 
\pr  Let $S:=B\setminus \M$, we have to check that $M_V$ is a model of $T_V$ and that the map sending $m\in M$ to $m+ann_S M$ is a morphism of $\L_V$-structures.
\qed
\medskip
\prop \label{gar} Let $M_V$ be a model of $T_{B,V}$.
Then we have an elementary embedding $M_{V}\hookrightarrow \prod_{\M} (M_{\M})_{V}$ as $\L_{V}$-structures.
\eprop
\pr By the pp elimination result for abelian structures, it is enough to show that given any pp $\L_{V}$-formula $\phi(\bar x)$ and $\bar a\in M$ such that $M_V\not\models \phi(\bar a)$, then for some maximal ideal $\M$ we have $(M_{\M})_{V}\not\models \phi(\bar a_{\M})$, where $\bar a_{\M}=\bar a+ann_{B\setminus \M}(M)$. 
\par Let $I:=\{r\in B:\;M_{V}\models \phi(\bar a.r)\}$. Then $I$ is a proper ideal of $B$; let $\M$ be a maximal ideal containing $I$ and $S:=B\setminus \M$. 
\par By the way of contradiction, suppose that $(M_{\M})_V\models \phi(\bar a_{\M})$. The formula $\phi(\bar x)$ is of the form $\exists \bar y\;\theta(\bar x, \bar y)$ where $\theta(\bar a,\bar y):=(\bar a.A_{1}+\bar y.A_{2}=0\;\&\;\bigwedge_{i} V_{\delta_{i}}(t_{i}(\bar a)+t'_{i}(\bar y)))$, with $\delta_{i}\in \Gamma^{+}$, and $A_{1},\;A_{2}$ are two matrices with coefficients in $B$. Let $\bar b=(d_{1}.s_{1}^{-1},\cdots,d_{n}.s_{n}^{-1})$, with $d_{i}\in M$, $s_{i}\in S$, $1\leq i\leq n$, be such that $(M_{\M})_{V}\models \theta(\bar a_{\M},\bar b_{\M})$. Equivalently, $\bar a.A_{1}+\bar b.A_{2}\in ann_{S} M$ and $\bigwedge_{i} V_{\delta_i}(t_{i}(\bar a_{\M})+t'_{i}(\bar b_{\M}))$. 
We multiply both expressions by $s=\prod_{i} s_{i}\in S$ and we get $\bar a.s.A_{1}+\bar b.s.A_{2}\in ann_{S} M$ and $\bigwedge_{i} V_{\delta_{i}.v(s)}(t_{i}(\bar a_{\M}.s)+t'_{i}(\bar b_{\M}.s))$ (using $(4)_V$). Since $\bar a.s.A_{1}+\bar b.s.A_{2}\in ann_{S} M$, there exists $\tilde s\in S$ such that $\bar a.s.\tilde s.A_{1}+\bar b.s.\tilde s.A_{2}=0$.
Finally we get: $\bar a.s.\tilde s.A_{1}+\bar b.s.\tilde s.A_{2}=0$ and $\bigwedge_{i} V_{\delta_{i}.v(s.\tilde s)}(t_{i}(\bar a_{\M}.s.\tilde s)+t'_{i}(\bar b_{\M}.s.\tilde s))$. Since $\delta_{i}.v(s.\tilde s)\geq \delta_{i}$, we get $V_{\delta_{i}}(t_{i}(\bar a.s.\tilde s)+t'_{i}(\bar b.s.\tilde s))$. Therefore, noting that $\bar b.s\in M$ and $M_{V}\models \theta(\bar a.s.\tilde s,\bar b.s.\tilde s)$, we obtain that $M_{V}\models \phi(\bar a.s.\tilde s)$. This shows that $s.\tilde s\in I\cap S$, a contradiction. \qed
\medskip
\subsection{Feferman-Vaught theorem}
\par Below, we introduce a property of the ring $B$ that implies the existence of relative complement for the basic closed sets in the Zariski spectrum of $B$.
\par Recall that for $c, d\in B$, we denoted  $c\in rad(d)$ the Jacobson radical relation, where $rad(d)$ is the intersection of all maximal ideals that contains $d$.
 \defn
Recall that $B$ has {\it good factorisation} \cite{DM}, if given any pair $a,\;b$ of non zeroes elements of $B$, there exist $c,\;d\in B$ such that $a=c.d$ with $gcd(c,b)=1$ and $b\in rad(d)$.
\edefn
\par First, let us link that last property with other possibly better known properties. 
\par If $B$ has good factorization, then given any two basic closed subsets $V(a), V(b)$ of $MSpec(B)$, there is an element $c$ such that  $V(c)=V(a)\setminus V(b)$. From that property it follows that the constructible
subsets of $MSpec(B)$ are either basic open or basic closed subsets \cite[Lemma 2.12]{DM}.  It also follows that if $B$ has good factorization, then $B$ is an elementary divisor ring \cite[Chapter III, Exercice 6.2]{FS}.

\defn A ring $R$ is adequate \cite[Exercice 6.4, page 118]{FS}, if for all $a\neq 0, b$, there exist $c,\;d$ such that $a=c.d$, $bR+cR=R$ and for all $d'$ such that $dR\subseteq d'R\subsetneq R$ implies $bR+d'R\subsetneq R$.
\edefn
\par It is easy to see that a B\'ezout ring $B$ with good factorization is adequate. Take $c, d\in B$ such that $a=c.d$ with $gcd(c,b)=1$ and $b\in rad(d)$ and let $d'$ be such that $dR\subset d'R$ and $d'R$ a proper ideal. Let $\M$ be a maximal ideal of $B$ containing $d'R$, so it contains $d$ and since $b\in rad(d)$, $b\in \M$. Therefore, $bR+d'R$ is a proper ideal of $B$.
\par Therefore, a B\'ezout ring with good factorization has the property that any prime ideal is contained in a unique maximal ideal \cite[Exercice 6.4, page 118]{FS}.
\medskip
\par Finally let us make the connection between {\it $B$ has good factorisation} and {\it $\Gamma(B)$ is a projectable $\ell$-group}. Let us first recall that last notion.
\par Given an $\ell$-group $\Gamma$, an ideal is a convex $\ell$-subgroup \cite[Section 3.2]{Gl} and a prime ideal is an ideal $P$ with the property that for any $f,\;g\in \Gamma$ such that $f\wedge g=1$, we have either $f\in P$ or $g\in P$ \cite[Section 3.3]{Gl}. By Zorn's Lemma, there exist minimal prime ideals. Let $Min(\Gamma)$ denotes the space of minimal prime ideals of $\Gamma$ endowed with the co-Zariski topology, namely the basis of open sets consists of $\{V(\delta):\;\delta\in \Gamma\}$ where $V(\delta):=\{P\in Min(\Gamma):\; \delta\in P\}$. 
\par Define for any $g\in \Gamma$, $g^{\perp}:=\{f\in \Gamma:\;\vert f\vert\wedge \vert g\vert=1\}$, where $\vert f\vert:=f\vee f^{-1}$ belongs to $\Gamma^{+}$. Recall that a {\it cardinal} sum of $\ell$-groups is a sum of $\ell$-groups endowed with the partial order defined componentwise \cite[Example 1.3.13]{Gl}. 
Then $\Gamma$ is {\it projectable} if for any $g\in \Gamma$, $g^{\perp}$ is a cardinal summand \cite[Section 3.5]{Gl}.
Since $MSpec(B)$ is homeomorphic to $Min(\Gamma(B))$  \cite[Proposition 8]{R}, the property for $B$ to have good factorization, translates into the property for $\Gamma(B)$ to be a projectable $\ell$-group. 
\bem\label{rad+} Using Remark \ref{rad} and the fact that $B$ is B\'ezout, one can further show that (\cite[Lemma 1.3 and its proof]{vdD}), for $a_{i}, b_{i}, c_{j}, d_{j}\in B$, $I$, $J$ finite, the following equivalences:
$$\forall \M\in MSpec(B)\;B_{\M}\models\;(\bigvee_{i\in I} a_{i}\nmid b_{i}\vee \bigvee_{j\in J} c_{j}\vert d_{j}),$$
$$\forall \M\in MSpec(B)\;(\bigwedge_{i\in I} (a_{i}:b_{i})\in \M\rightarrow \prod_{j\in J} (c_{j}:d_{j})\in \M),$$
$$B\models\;\prod_{j\in J} (c_{j}:d_{j})\in rad(gcd((a_{i}:b_{i})_{i\in I})),$$
where $gcd((a_{i}:b_{i})_{i\in I})$ denotes a generator of the ideal generated by the elements $(a_{i}:b_{i})$, $i\in I.$
\ebem
\bigskip
\thm \label{ppelB} 
Let $M_{V}$ be a model of $T_{B,V,div}$. Let $\phi(\bold{y})$ be a pp $\L_{V,B}$-formula. Then there are finitely many conjunctions of atomic (c.a.) $\L_{V,B_{\M}}$-formulas $\psi_{k}(\bold{y})$, $k\in K$, and a finite covering of $MSpec(B)$ into constructible subsets $C_{\phi,k}$
such that for any $\bold{u}\in M$, we have: 
$$M_{V}\models \phi(\bold{u}) \leftrightarrow (\bigwedge_{k\in K}\;{\rm for\;all}\; \M\in C_{\phi,k}\;(\M_{\M})_{V}\models \psi_{k}(\bold{u}_{\M})),$$
where $\bold{u}_{\M}$ denotes the image of the tuple $\bold{u}$ in $M_{\M}$.
\ethm
\pr First, by Proposition \ref{gar}, we have $M_{V}\models \phi(\bold{u})$ iff for all $\M\in MSpec(B)$, $(M_{\M})_{V}\models \phi(\bold{u}_{\M})$. 
Each $M_{\M}$ is a $\L_{B_{\M},V}$-structure and since $B_{\M}$ is a valuation domain, we may apply Theorem \ref{ppel} to these classes of $\L_{B_{\M},V}$-structures $M_{\M}$. 
So there exist finitely many constructible subsets $C_{\phi,k}$, $k\in K$, with $K$ finite such that for any $\M\in C_{\phi,k}$, $Mod_{B_{\M}}\models \forall \bold{y}(\phi(\bold{y})\leftrightarrow \psi_{k}(\bold{y}))$, where $\psi_{k}$ is a  c.a. $\L_{V, B_{\M}}$-formula. \qed
\medskip
\par In case $B$ has good factorization, we obtain a neater statement.
\cor \label{cons} Assume that $B$ has good factorization. 
Let $M_{V}$ be a model of $T_{B,V,div}$. Let $\phi(\bold{y})$ be a pp $\L_{V,B}$-formula. Then there are finitely many conjunctions of atomic $\L_{V,B_{\M}}$-formulas $\theta_{i}(\bold{y})$, $i\in I$, and a finite partition of $MSpec(B)$ into basic open or basic closed subsets $O_{i}$, such that for any $\bold{u}\in M$, we have: 
$$M_{V}\models \phi(\bold{u}) \leftrightarrow (\bigwedge_{i\in I}\;{\rm for\;all}\; \M\in O_{i}\;(\M_{\M})_{V}\models \theta_{i}(\bold{u}_{\M})),$$
where $\bold{u}_{\M}$ denotes the image of the tuple $\bold{u}$ in $M_{\M}$.
\ecor
 \pr Since $B$ has good factorisation, any constructible subset of $MSpec(B)$ is either a basic open or basic closed subset of $MSpec(B)$ \cite[Lemma 2.12]{DM}.\qed
 \medskip
 \prop \label{dec} Suppose that a $B$ is a countable B\'ezout domain 
 and assume that for each $\M\in MSpec(B)$, the quotient $B/\M$ is infinite. 
 Further, assume that $B$ satisfy hypothesis $(EF)$ and that the Jacobson radical relation $rad$ is recursive. Then $T_{B}$ is decidable.
\eprop
\pr First recall that the theory $T_{B,V,div}$ is a definable expansion by definition of the theory $T_B$. The key ingredient is Theorem \ref{ppelB} which, given a pp $\L_V$-formula $\phi(x)$, enables us to obtain (in an effective way) a finite covering of $MSpec(B)$ into constructible subsets: $C_{\phi,k}$, $k\in K$ and finitely many c.a. $\L_{V,B_{\M}}$-formulas $\psi_{k}(x)$ such that over each $C_{\phi,k}$, $\phi(x)$ 
  is equivalent to $\psi_{k}$ which can be assumed to be of the form: $x.a_{k}=0\;\&\;\bigwedge_{i\in I}\;V_{\delta_{k_i}}(x.b_{k_i})$, $a_{k},\;b_{k_i}\in B_{\M}, \delta_{k_i}\in \Gamma_{\M}^+$ by Theorem \ref{ppel} $(\star)$.
\par Then we use a standard procedure to obtain decidability of the theory $T_{B,V,div}$ (see for instance \cite[Theorem 6.2]{PPT}), that we detail below.
\par The hypothesis on the ring $B$ implies that the theory $T_{B,V,div}$ is recursively enumerable. 
As recalled in Remark \ref{Zg}, since for every maximal ideal $\M$, $B/\M$, is infinite, proving that $T_{B,V,div}$ is decidable is equivalent to being able to answer the question whether in $Zg_B$, a basic open set is included in a given finite union of other basic open sets, namely $[\phi_{0}/\psi_{0}]\subset \bigcup_{i=1}^n [\phi_{i}/\psi_{i}]$, with $\phi_{i},\psi_{i}$ pp $\L_{V}$-formulas, and $\psi_{i}\rightarrow \phi_{i}$, $0\leq i\leq n$  $(\star\star)$. 
A point in the Ziegler spectrum is (the isomorphism class of) an indecomposable pure-injective $B$-module and so a $B_{\M}$-module for some maximal ideal $\M$ of $B$ \cite[Theorem 5.4]{Z}. 
\par Given the above finite set of pp $\L_{V,B}$-formulas $\phi_{i},\psi_{i}$, $0\leq i\leq n$, we apply to each of these formulas procedure $(\star)$ and we obtain (in an effective way) a finite covering of $MSpec(B)$ into constructible subsets: $C_{\ell}$, $\ell\in L$, over which each of these pp $\L_{V,B}$-formulas is equivalent to a c.a. formula of the form: $x.a_{\ell}=0\;\&\;\bigwedge_{i\in I}\;V_{\delta_{\ell_i}}(x.b_{\ell_i})$, $a_{\ell},\;b_{\ell_i}\in B_{\M}, \delta_{\ell_i}\in \Gamma_{\M}^+$. 
Moreover using duality, by Corollary \ref{dec-val} (see Claim \ref{1}), we may only consider open sets of the form $[x.b=0\,\&\;V_{\delta_1}(x)/V_{\delta_2}(x.c)]$, $b, c\in B_\M$, $\delta_1, \delta_2\in \Gamma_\M^+$.
\par In order to check whether $(\star\star)$ holds, we proceed then as in the proof of Corollary \ref{dec-val} (see Claim \ref{5}), and it reduces on each element $C_{\ell}$, $\ell\in L$, of the covering of $MSpec(B)$, to divisibility conditions on elements of $B$ in the localizations and order relations between the $\delta$'s, which reduce to divisibility conditions in the corresponding $B_{\M}$. 
\par Finally, using Remark \ref{rad+}, this can be expressed using the Jacobson radical relation in $B$. We have to answer statements which are finite conjunctions of the following form: 
$\prod_{i\in I'} (r'_i:s'_i)\in rad (gcd(a'_j:c'_j)_{j\in J'})$, where the (finite) index sets $I',\;J'$ and the elements $a'_{j}, c'_{j}, r'_{i}, s'_{i}$ can be effectively determined from the previous data. Then we use the hypothesis $(EF)$ on our ring to effectively obtain the elements $(r'_i:s'_i),\;(a'_j:c'_j)$, $gcd(a'_j:c'_j)_{j\in J'}$ from the previous ones. Finally we use the hypothesis that the Jacobson radical relation is recursive in order to decide whether $\prod_{i\in I'} (r'_i:s'_i)$ belongs to $rad (gcd(a'_j:c'_j)_{j\in J'})$.
\qed
\section{Applications}\label{App}
In this section we will revisit in details the examples mentioned in the introduction: on one hand the ring of algebraic integers and on the other hand the ring of holomorphic functions over $\IC$ and we also examine the cases of real and $p$-adic algebraic integers. 
\subsection{Good Rumely domains} 
In order to axiomatize the elementary theory of the ring $\widetilde \Z$ of algebraic integers, the following subclasses of B\'ezout domains were introduced in \cite{DM}.
A domain $B$ with fraction field $K$ is a Rumely domain if it has the following properties:
\begin{itemize}
\item[(Ru 1.)] The field $K$ is algebraically closed.
\item[(Ru 2.)] Every finitely generated ideal of $B$ is principal.
\item[(Ru 3.)] (Local-global principle) If $C \subseteq  \mathbb{A}^m(K)$ is a smooth irreducible closed curve, $f \in  K[X_1, \ldots , X_m]$ and $C_f = \{x \in C : f(x)\neq 0 \}$ has points in $(1/a)\mathbb{A}^m(B)$ and in $ (1/b)\mathbb{A}^m(B)$, where $a, b \in B\setminus \{ 0 \}$ are relatively prime, then $C_f$ has a point in $\mathbb{A}^m(B)$.
\end{itemize}
$B$ is a good Rumely domain if it satisfies, moreover, the following properties.
\begin{itemize}
\item[(Ru 4.)]  (Good factorization) For all $a,\;b\in B^{\star}$, there are $c,\;a_{1}\in B$ such that $a=c.a_{1}$ with $gcd(c,b)=1$ and $b\in rad(a_{1})$.
\item[(Ru 5.)] Every nonzero nonunit is the product of two relatively prime nonunits.
\item[(Ru 6.)] Its Jacobson radical, namely the intersection of all maximal ideals of $B$, is equal to $\{ 0 \}$ and $B \neq K$.
\end{itemize}
All these properties are first-order expressible in the language of rings (\cite[1.6]{DM}). 
Consider the Boolean algebra $\B(B)$ generated by all basic closed subsets $V(a)$, $a\in B$, of the maximal spectrum $MSpec(B)$; if a B\'ezout domain $B \neq K$ satisfies (Ru 4), then one can check that
(Ru 5) holds in $B$ if and only if $\B(B)$ is an atomless boolean algebra.
\par In \cite{PS}, A. Prestel and J. Schmid axiomatize a class of (commutative) domains endowed {\it a radical relation $\preceq$} \cite[Introduction]{PS}. Alternatively they show that for each such relation, one can associate a non-empty subset $P_{\preceq}$ of the prime spectrum in such a way $a\preceq b$ if and only if for every prime ideal $I\in P_{\preceq}$, $a\in  I\;\rightarrow b\in I$ \cite[Theorem 2.5]{PS}. They show that the class of good Rumely domains is exactly the class of existentially closed (e.c.) $r_0$-domains $(R,\preceq)$, namely those $(R,\preceq)$ such that $(0\preceq a\rightarrow a=0)$ \cite[Theorem 3.3]{PS}.
They also note that in an e.c. $r_0$-domain $R$, the relation $\preceq$ is induced by the maximal spectrum $MSpec(R)$ of $R$, namely $a\preceq b$ if and only if  $V(a)\subseteq V(b)$. 
\medskip
\par Examples of good Rumely domains are: the ring $\widetilde \Z$ of algebraic integers, the integral closure of $\F_{p}[t]$.

\par All localizations of Rumely domains are again Rumely domains (\cite[Corollary 3.5]{DM}). Localizations of B\'ezout domains with good factorisation are again B\'ezout domains with good factorisation (\cite[2.10]{DM}). 
\par Let $\O$ be either the ring of algebraic integers in a number field or the integral closure of $\F_p[t]$ in a finite degree field extension of $\F_p(t)$ and let $S$ denote a multiplicative subset of $\O$. Then assume that the ring $S^{-1}.\O$ is not a field and that $S^{-1}.\O$ has infinitely many maximal ideals, then the Jacobson radical of $\widetilde{(S^{-1}.\O)}$ is $\{0\}$ \cite[Lemma 2.15]{DM}. 
\par Note that when $B$ is a good Rumely domain, the prime radical relation and the Jacobson radical relation coincide \cite[Theorem 3.3]{PS} and L. Gregory showed that the decidability of the theory of $B$-modules
implies that the prime radical relation is recursive \cite[Lemma 3.2]{Gr}. Therefore we get the following Corollary to Proposition \ref{dec}.
\cor Let $B$ be a countable good Rumely domain, assume that 
 $B$ satisfies (EF).  Then $T_{B}$ is decidable if and only if the prime radical relation is recursive.  \qed
\ecor 
\pr The only thing to note is that  for each $\M\in MSpec(B)$, the quotient $B/\M$ is infinite. If the characteristic of $B$ is zero, this is immediate and if the characteristic of $B$ is a prime $p$, all good Rumely domains containing the prime field $\F_p$ also contains its algebraic closure \cite[Theorem 4.2]{DM}.
\qed
\par The field $\widetilde \Q$, the algebraic closure of $\Q$, can be equipped with a recursive structure \cite[page 131]{Ru} and from that presentation one can deduce that the ring $\widetilde \Z$ can also be equipped with a recursive structure. Earlier, M. Rabin showed that if $F$ is a computable field, then so is its algebraic closure \cite[Theorem 7]{Rabin}.

\bem \label{comp}{\rm  \cite[Fact 2]{vdD}} Suppose the ring $R$ satisfies hypothesis (EF) and that the Jacobson radical relation is equal to the prime radical relation, then the Jacobson radical relation $rad$ on $R$ is recursive.
\ebem
\pr For the reader convenience, we give the proof below \cite[page 192]{vdD} (note that van den Dries uses that the ring is equipped with a recursive structure, but in our context, we may replace this by assumption (EF)). Van den Dries uses that the relation $x\in rad(y_1,\cdots,y_{\ell})$ is recursively enumerable (r.e.) as well as its complement.
To show it is r.e., one writes: $x\in rad(y_1,\cdots,y_{\ell})\leftrightarrow \exists n\in \N\;x^n\in (y_1,\cdots,y_{\ell})$ and that its complement is r.e.
$x\notin rad(y_1,\cdots,y_{\ell})\leftrightarrow \exists z\;(1\in (z,x)\;\&\;1\notin(z,y_1,\cdots,y_{\ell})).$ Since our ring is B\'ezout, we have that $1\in (z,x)\leftrightarrow gcd(z,x)=1$.
\qed
\par Therefore,  one can deduce the following Corollary.
\begin{corollary} The theory
 $T_{\widetilde \Z}$ is decidable. \qed
\end{corollary}
 \par The decidability of the theory of modules over the ring of algebraic integers $\widetilde \Z$  has also been obtained by S. L'Innocente, G. Puninski and C. Toffalori, using different methods \cite{LPT}.
 \subsection{Real algebraic integers and p-adic integers} 
A. Prestel and J. Schmid used the same analysis as described above (for $\widetilde \Z$) in order to study the rings $\widetilde \Z\cap \IR$ and $\widetilde \Z\cap \Q_{p}$ (\cite{PS91}). They showed that in the case of $\widetilde \Z\cap \IR$ and $\widetilde \Z\cap \Q_{p}$, the axiomatizability depends on a certain local-global principle (as in the case of $\widetilde \Z$). Furthermore in these two rings, any prime ideal is maximal since it holds in $\Z$ \cite[Corollary 5.8, page 61]{AM}. By working in the setting of rings $(R, \preceq)$ with a radical relation $\preceq$, they proved that the related theories of  $\widetilde \Z\cap \IR$ and $\widetilde \Z\cap \Q_{p}$ (in the language of rings) are decidable \cite[Corollary 2.5 and Corollary 3.5]{PS91}. So, Proposition \ref{dec} leads us to the following corollary.

\begin{corollary}
Let $B$ be one of the two rings $\widetilde \Z\cap \IR$ and $\widetilde \Z\cap \Q_{p}$.
Then the corresponding theory $T_{B}$ is decidable. \qed
\end{corollary}
\pr It suffices to show that each of these rings satisfy hypothesis (EF), by Remark \ref{comp}. Moreover since $\widetilde \Z$ satisfies hypothesis (EF) \cite[Fact 2]{vdD}, it remains to check that the intersections $\widetilde \Z\cap \IR$ and $\widetilde \Z\cap \Q_{p}$ are recursive. As in \cite{vdD}, we will use the recursive structure on $\widetilde \Z$ defined by Rumely \cite[III, page 131]{Ru}.
\par As in \cite[III, page 131]{Ru}, we fix an embedding of $\widetilde \Q$ in $\IC$ and represent each element $\alpha$ of $\widetilde \Z$ as a pair $(f_{\alpha}(x), a+bi)$ where $f_{\alpha}(x)\in \Z[x]$ is the minimal monic polynomial of $\alpha$ and $a+b.i\in \IC$ is a {\it sufficiently good} decimal approximation of $\alpha$ to distinguish it from its conjugates. (There is a discussion in \cite[page 132]{Ru} to how {\it good} is an approximation good enough.) 
\par First consider $\widetilde \Z\cap \IR$. One can give an estimate of the minimal distance $B(f_{\alpha})$ of each of the roots of $f_{\alpha}(x)$ (in terms of the coefficients of $f_{\alpha}(x)$) and in order to check that $\alpha\in \widetilde \Z\cap \IR$, we express that the complex conjugate $a-bi$ is also a root of $f_{\alpha}(x)$ at distance strictly smaller than $B(f_{\alpha})$.
\par In case of $\widetilde \Z\cap \Q_p$, we can proceed as follows. By the result of A. Prestel and J. Schmid recalled above, the theory of the ring $\widetilde \Z\cap \Q_{p}$ is decidable.  
So, we can check whether the sentence $\exists x\,(f(x)=0\;\&\;\vert x-(a+bi)\vert <B(f_{\alpha}))$ holds in $\widetilde \Z\cap \Q_{p}$. If the answer is {\it yes}, we keep such $\alpha$.
\qed
\subsection{The ring of entire functions and its integral closure}\label{entire}

\par Let $B$ be the ring $\H(\IC)$ of entire functions in $\IC$.  As we already recalled (see Example \ref{hol}), $\H(\IC)$ is a B\'ezout domain, as is its integral closure $\widetilde \H(\IC)$. Moreover $\widetilde \H(\IC)$ satisfy the algebraic properties $(Ru 1.)$ up to $(Ru 6.)$ listed above, except the property $(Ru 3.)$ of the local-global principle \cite[5.6]{DM}.  Here we will restrict ourselves to $\H(\IC)$. 
 \lmm The B\'ezout domain $\H(\IC)$ has good factorization.
 \elmm
\pr By Weierstrass factorization theorem \cite[Theorem 15.10]{Rudin}, one can write $f$ as $e^{h}.z^d.\prod_{n=1}^{\infty} E_{n-1}(\frac{z}{z_{n}})$, where $Z(f)\setminus\{0\}=\{z_n:\;n\in \omega\}$ and $d$ is the multiplicity of $0$ as a zero of $f$. Then let $Z_1=\{n\in \omega:\;z_n\in Z(f)\cap Z(g)\}$ and set, if $g(0)\neq 0$, 
$f_1:=e^{h}.\prod_{n\in Z_1} E_{n-1}(\frac{z}{z_{n}})$,  $f_2:=z^d.\prod_{z_n\notin Z(g)} E_{n-1}(\frac{z}{z_{n}})$ and if $g(0)=0$, set 
$f_1:=e^{h}.z^d.\prod_{n\in Z_1} E_{n-1}(\frac{z}{z_{n}}),\; f_2:=\prod_{z_n\notin Z(g)} E_{n-1}(\frac{z}{z_{n}})$. We have $f=f_1.f_2$ and if  $\M$ is any maximal ideal containing $f_1$, 
since $Z(f_1)\subset Z(g)$, we get $g\in \M$ by $(\star)$.\qed
\medskip
\par I. Kaplansky noted that there are prime non-maximal ideals in $\H(\IC)$ \cite[Theorem 1]{H}. A necessary and sufficient condition for a prime ideal $P$ to be non-maximal is that for all $f\in P$, the multiplicity function $\mu_f$ is unbounded ($\mu_f$ as in Example \ref{hol}).
\medskip
 \par Since $\H(\IC)$ is uncountable, there is the usual problem of defining a suitable notion of decidability of a theory of modules in that case. One could take $R$  a countable elementary substructure of $\H(\IC)$ (in the language of rings) (respectively of $\widetilde \H(\IC)$) and assume that $R$ is effectively given and that the Jacobson radical relation is recursive. 
From Proposition \ref{dec}, we get that the corresponding theory $T_{R,V}$ is decidable. Of course it would be more informative to exhibit such a subring. In a forthcoming paper with G. Puninski and C. Toffalori, we describe the Ziegler spectrum over $\H(\IC)$ \cite{LPPT}.
\par Finally let us note that, contrary to the other examples of rings we considered, it is an open question whether the positive existential theory of $\H(\IC)$ in the language of rings expanded with a new constant symbol interpreted by the identity function of $\H(\IC)$ is decidable \cite[Problem 1.1]{GP}.
\bigskip
\section{Further work}
\par Now, we introduce the notion of $\ell$-valued $B$-modules which extends the notion of valued modules occurring in, for instance \cite{Co},  \cite{F} or \cite{L}(\S2) and also in \cite{Dries}, \cite{BP}, for a model-theoretic point of view.  
Let $M$ be a $B$-module and set $M^{\star}:=M\setminus \{0\}$. Let $\Gamma:=\Gamma(B)$ be the $\ell$-group of divisibility of $B$ (with the group law $\cdot$, lattice operation $\wedge$ and neutral element $1$) (see section 2.1).
\defn \label{W} 
A  {\em $\ell$-valued $B$-module} is a two-sorted structure 
$(M,\bar \Gamma(B)^+,w)$, where 
$M$ is a $B$-module 
and $w: M\to \bar \Gamma(B)^+$ such that 
\begin{enumerate}
\item for all $m_1, m_2\in M$, $w(m_1+m_2)\ge
w(m_1)\wedge w(m_2)$, and $w(0)=\infty$; 
\item for all $m\in M^{\star}$, $w(m. a)=w(m)\cdot v(a)$, for each $a\in B^{\star}$.
\end{enumerate}
\edefn
\bsp Considering $B$ as a module over itself and letting $v$ the $\ell$-valuation on its group of divisibility, we get that $B$ is a
$\ell$-valued $B$-module.
\ebsp
\bem We could have taken $\Gamma$ any $\ell$-group, or even we could have considered a distributive lattice $\Delta$, assuming that for each of them we have an action of $\Gamma(B)$, as in, for instance, \cite{BP}. 
\ebem
\medskip
\par From the axioms above, we easily deduce the 
following properties:
\par $(P.1)$ $w(m)=w(-m)$ and for all $m\in M$, $w(m)\leq w(m)\cdot s$, for each $s\in \Gamma(B)^+$  \label{increasing}.
\par $(P.2)$ Let $m,\;n\in M$ and assume that $w(n)\geq w(m)$. We have $w(m+n)\geq w(m)\wedge w(n)=w(m)$ and $w(m)\geq w(m+n)\wedge w(-n)=w(m+n)\wedge w(n)=w(m+n)$. So, $w(m+n)=w(m)$.
\medskip

\par Given $(M,\bar\Gamma,w)$ an  $\ell$-valued $B$-module; we may consider it in the weaker formalism of abelian structures. Namely we associated with it  
the $\L_{V}$-structure $M_{V}$ where $M_{\gamma}:=\{m\in M : w(m)\ge \gamma\}$, $\gamma\in \Gamma^+$. It is easily checked that $M_V$ is a model of $T_{B,V}$. 
\bigskip
\par {\bf Acknowledgments:} The second author would like to thank the Mathematical Department of Camerino University and in particular Carlo Toffalori, for their hospitality in the fall 2015 during which the present work was done. We would like to thank D. Macpherson and G. Puninski for their advice on the writing up of these results.  

\end{document}